 \documentclass[authoryear,preprint,review,12pt]{elsarticle}

\usepackage{amssymb}
\usepackage{lipsum}
\usepackage{color}
\usepackage{amsmath,bm}
\usepackage{optidef}
\usepackage{cleveref}

\usepackage{algorithm}
\usepackage{algorithmic}
\usepackage{graphicx}
\usepackage{multirow}
\usepackage{subcaption}
\usepackage{pgfplots}
\usepackage{tikz}
\usetikzlibrary{calc}
\usepackage{amsmath} 
\usepackage{lscape}
\usetikzlibrary{spy}
\usepackage{geometry}
\usepackage{tabularx}
\usepackage{adjustbox}
\usepackage{xcolor}
\usepackage{natbib}
\definecolor{deepbrown}{rgb}{0.36, 0.25, 0.20}

\journal{}

\pgfplotsset{compat=1.18} 
\begin{document}

\begin{frontmatter}

\title{Multi-period railway line planning for integrated passenger-freight transportation\tnoteref{title}}

\author[first]{Wanru Chen}
\author[second]{Rolf N. van Lieshout}
\author[first]{Dezhi Zhang\corref{cor1}}
\ead{dzzhang@csu.edu.cn}
\author[second]{Tom Van Woensel}
\affiliation[first]{{School of Traffic and Transportation Engineering, Central South University, }, 
            city={Changsha},
            state={Hunan},
            postcode={410075}, 
            country={China}}
\affiliation[second]{{Department of Operations, Planning, Accounting, and Control, School of Industrial Engineering, Eindhoven University of Technology, },
            city={Eindhoven},
            state={North Brabant},
            postcode={5600 MB}, 
            country={Netherlands}}
\cortext[cor1]{Corresponding author.}

\begin{abstract}
This paper addresses a multi-period line planning problem in an integrated passenger-freight railway system, aiming to maximize profit while serving passengers and freight using a combination of dedicated passenger trains, dedicated freight trains, and mixed trains. To accommodate demand with different time sensitivities, we develop a period-extended change\&go-network that tracks the paths taken by passengers and freight. The problem is formulated as a path-based mixed integer programming model, with the linear relaxation solved using column generation. Paths for passengers and freight are dynamically generated by solving pricing problems defined as elementary shortest-path problems with duration constraints. We propose two heuristic approaches: price-and-branch and a diving heuristic, with acceleration strategies, to find integer feasible solutions efficiently. Computational experiments on the Chinese high-speed railway network demonstrate that the diving heuristic outperforms the price-and-branch heuristic in both computational time and solution quality. Additionally, the experiments highlight the benefits of integrating freight, the advantages of multi-period line planning, and the impact of different demand patterns on line operations.
\end{abstract}




\end{frontmatter}




\section{Introduction}
With the continuous development of the global economy and e-commerce, freight transportation has rapidly increased, especially for express delivery services. This surge places significant pressure on road networks and contributes to excessive environmental pollution, highlighting the need for transportation solutions that are cost-effective, environmentally friendly, and reliable. One promising strategy to address this issue is the integration of freight with railway passenger transportation, which has been successfully implemented by, for example, CRH Express in China, La Poste in France, Parcel InterCity in Germany, and Mercitalia Fast in Italy.

\begin{figure}[t]
    \centering
    \begin{tikzpicture}[scale=0.8, every node/.style={scale=0.6}, >=stealth]
    \draw (0, 0) rectangle (2, 1);
     \node at (1, 0.5) {\parbox{1.8cm}{\centering\baselineskip=12pt Network\\Design}};
    \draw[->,thick] (2,0.5) -- (2.5,0.5);
    \draw (2.5, 0) rectangle (4.5, 1);
    \node at (3.5, 0.5) {\parbox{1.8cm}{\centering\baselineskip=12pt Line\\Planning}};
    \draw[->,thick] (4.5,0.5) -- (5,0.5);
    \draw (5, 0) rectangle (7, 1);
    \node at (5.95, 0.5) {\parbox{2cm}{\centering Timetabling}};
     \draw[->,thick] (7,0.5) -- (7.5,0.5);
     \draw (7.5, 0) rectangle (9.5, 1);
    \node at (8.5, 0.5) {\parbox{1.8cm}{\centering\baselineskip=12pt Traffic \\Planning}};
    \draw[->,thick] (9.5,0.5) -- (10,0.5);
     \draw (10, 0) rectangle (12, 1);
    \node at (11, 0.5) {\parbox{1.8cm}{\centering\baselineskip=12pt Rolling \\Stock}};
     \draw[->,thick] (12,0.5) -- (12.5,0.5);
     \draw (12.5, 0) rectangle (14.5, 1);
    \node at (13.5, 0.5) {\parbox{1.8cm}{\centering\baselineskip=12pt Crew\\Scheduling}};
    \draw[->,thick] (14.5,0.5) -- (15,0.5);
    \draw (15, 0) rectangle (17, 1);
    \node at (16, 0.5) {\parbox{1.8cm}{\centering\baselineskip=12pt Shunting}};
    \draw[->,thick] (17,0.5) -- (17.5,0.5);
     \draw (17.5, 0) rectangle (19.5, 1);
    \node at (18.5, 0.5) {\parbox{1.9cm}{\centering\baselineskip=12pt Real-Time\\Control}};
   
    \node at (1.25, -0.8) {\large Strategic Level};
    \draw[->,thick] (2.75,-0.8) -- (8.8,-0.8);
    \node at (10.25, -0.8) {\large Tactical Level};
    \draw[->,thick] (11.75,-0.8) -- (16.25,-0.8);
    \node at (18, -0.8) {\large Operational Level};
\end{tikzpicture}
        \caption{The different planning stages of a railway operator}
    \label{F_LP}
\end{figure}
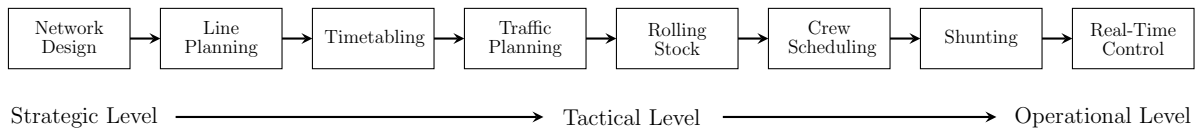

Integrating freight in passenger railway transportation requires adjustments to the planning process railway operators face, which spans from the strategic to the operational level, as shown in Figure \ref{F_LP}. Existing approaches for freight integration mainly start from tactical adaptations and include the optimization of timetables, rolling stock, and shunting schedules \citep{R8,R21,R14,R7,R15,R6,R18}. However, these adjustments cannot fully exploit the opportunities that combining passengers and freight in one system offers, as this systemic change also requires consideration in strategic planning. To address this gap, we focus on the line planning problem, which establishes the routes and frequencies for trains and serves as input for the timetabling process. Solving this fundamental strategic issue can optimize resource and infrastructure use to meet the long-term demand of combined transportation systems, enhancing efficiency and maximizing profitability.

When considering line planning for an integrated passenger-freight system, it is crucial to take into account different time sensitivities of passenger and freight demand. Passenger demand is highly time-sensitive, with strong preferences for specific departure or arrival times and short travel times. This leads to significant demand variations during the day, with traffic flow peaks during morning and evening rush hours \citep{L12}. Conversely, freight demand is less time-sensitive and can be scheduled flexibly throughout the day and with longer travel times. Most existing line planning literature assumes that the same line plan is operated throughout the day, which cannot sufficiently accommodate these different time-sensitivities. Therefore, we propose a multi-period line planning problem that allows us to efficiently manage time-varying passenger demand while using available spare capacity in different periods to accommodate freight demand. 

Additionally, freight integration requires accurate assignment of passengers and freight to lines to ensure that offered capacity is used efficiently but never exceeded. This is especially challenging in line planning since, on the one hand, passengers and freight rely on established lines to determine their paths. On the other hand, optimizing the lines requires knowing the passengers and freight paths. To address this, we develop a so-called period-extended change\&go-network that traces the paths taken by passengers and freight through space (on the line level) and time (on the period level). Service quality is ensured by constraining the durations of the paths taken by passengers and freight compared to the shortest possible paths. 

We formulate the multi-period line planning problem for an integrated passenger and freight system as a path-based mixed integer program (MIP). This formulation is hard to solve directly for practical instances due to the many potential paths involved. Hence, column generation is employed to solve its linear programming (LP) relaxation, and two heuristics, namely price-and-branch and a diving heuristic, are proposed to find integer solutions efficiently. Furthermore, we develop acceleration strategies to reduce the computation time of these heuristics.

We test our approach using instances derived from a Chinese high-speed railway network comprising 42 stations. Multiple operational scenarios are considered, each involving a combination of dedicated passenger trains, dedicated freight trains, and trains transporting passengers and freight. Across all instances, the diving heuristic achieves a maximum average gap of 1.6\% in smaller instances and 2.8\% in larger instances, consistently outperforming the price-and-branch heuristic. Management insights highlight that integrating freight into a multi-period line plan can enhance profitability while maintaining or even increasing passenger service levels. However, different demand patterns can affect the overall service capacity of the system, as well as the capacity allocation between passenger and freight operations.

To summarize, the contribution of our paper is threefold. First, we present a multi-period line planning problem and demand assignment approach that considers different time sensitivities of passenger and freight demand, enhancing profitability through long-term freight integration planning. Second, we introduce column generation-based heuristics to solve the problem, leveraging the structure of the period-extended change\&go-network and using acceleration strategies to speed up the solution method. Third, we deliver critical managerial insights for railway operators through comprehensive testing using practical instances. We illustrate the algorithm's performance in several scenarios and present optimization results that confirm the effectiveness of our method in enhancing decision-making processes.

The remainder of this paper is organized as follows: Section 2 reviews the literature on integrating passenger and freight in railway scheduling and the line planning problem. Section 3 outlines the problem's features and the proposed optimization model. Section 4 introduces two column generation-based heuristics. We present the experimental design in Section 5, followed by the computational performance and management insights in Sections 6 and 7. Finally, Section 8 presents the conclusions and future outlook of the study.

\section{Literature review}
In the following, we discuss the literature on the multi-period railway line planning problem for integrated passenger-freight transportation. First, we review existing works on integrating passenger and freight in railway scheduling. Second, we present papers that specifically address line planning.

\subsection{Integrating passenger and freight in railway scheduling}
Integrating passenger and freight services in railway scheduling represents a significant shift from traditional practices, where these operations are managed separately. Integration has the potential to not only optimize operations and increase profits but also deliver environmental benefits \citep{R2,R17,R19}. Innovative and collaborative approaches at strategic, tactical, and operational levels are essential to realize these benefits.

At the strategic level, existing research primarily emphasizes long-term planning involving train routing, capacity allocation, and network design. For instance, \cite{R13} examines the impact of the North Doncaster Chord on freight operations in the British railway network, mainly focusing on train routing, scheduled journey times, and train punctuality. \cite{R3} present a two-level passenger and freight train planning model for shared-use railway corridors, considering passenger access priority. \cite{R9} develop an analytical model to optimize track allocation, consolidation time, and pricing levels for shared-use railway corridors serving both passenger and freight trains, considering stochastic demand and distinguishing between train characteristics. A two-stage stochastic programming model is proposed to address a network planning and freight assignment problem by \cite{R20}  for high-speed railway express systems, considering intermodal station selection, vehicle assignment, and freight flow organization. 

At the tactical and operational level, the focus shifts to medium-term or short-term decisions involving timetable planning and resource allocation, aiming to optimize the utilization of available infrastructure while minimizing conflicts between passenger and freight transportation. Both \cite{R7} and \cite{R16} address the train capacity allocation problem in railway systems for mixed passenger and freight transportation, with \cite{R7} specifically considering uncertain demand. For urban freight transportation, utilizing passenger railway networks as environmentally friendly alternatives is proposed, with studies primarily addressing the optimization of freight flow distribution within service networks originally designated for passenger traffic. \citep{R14,R15,R6,R21,R18}, while in the research by \cite{R4} and \cite{R5}, freight can be transported by utilizing the extra space inside the passenger trains or inserting dedicated freight trains. Furthermore, \cite{R8} explore the insertion of freight trains into existing timetables to achieve a harmonious coexistence with passenger transportation in railway networks.  

Strategic planning plays a crucial role in establishing the foundation for successful freight integration within passenger transport. However, there remains a significant gap in research concerning line planning. Addressing this gap is crucial for developing comprehensive solutions that effectively coordinate passenger and freight operations, considering their unique features.

\subsection{Line planning problem}
Two main orientations for line planning problems are service-oriented and cost/profit-oriented. Service-oriented studies prioritize providing more direct travels or reducing travel/riding times \citep{L18,L21, L8}. Cost/profit-oriented studies aim to maximize profit or minimize total costs \citep{L5,L6,L3,L2,L9,L16}. Additionally, some research simultaneously considers these optimization objectives, balancing efficiency and profitability while meeting passenger needs \citep{L4, L20, L12}.

Most studies assume a predetermined set of allowed lines, called a \textit{line pool}. Some integrate line route generation directly into decision-making processes. For example, \cite{L7} introduce models allowing lines to have varying halting patterns instead of stopping at all stations. Additionally, \cite{L4} propose a multi-commodity flow model enabling dynamic line generation. In contrast, \cite{L1} introduce a bi-level multi-objective mixed-integer nonlinear programming model to optimize passenger assignment and train routing simultaneously in high-speed railway line plans.

In line planning, demand is typically given in matrix form for origin-destination (OD) pairs, necessitating decisions on how this demand is routed throughout the network. In the literature, demand assignment is approached in two main ways. The first method defines paths independently of lines, such as on the physical transportation network. This method ensures that all links in the network can handle the distributed demand with sufficient service capacity; however, it does not specify how demand is assigned to lines \citep{L18, L4, L1}. In contrast, the second method defines paths based on lines, such as on the change\&go-network. This allows direct allocation of demand to specific lines, making it possible to track the lines taken and the exact transfers \citep{L9, L20, L8}. The first method ensures overall demand coverage broadly, while the second method offers detailed demand assignment, which helps in precise line and resource optimization.

\begin{table}[t]
\renewcommand{\arraystretch}{0.9}
\centering
\caption{Comparison of key elements across different studies on line planning problem}
\adjustbox{max width=1\textwidth}{
\begin{tabular}{lllcccccccccl}
\hline
\multirow{2}{*}{References} &  & \multirow{2}{*}{\begin{tabular}[c]{@{}l@{}}Objective \\ orientation\end{tabular}} &  & \multirow{2}{*}{\begin{tabular}[c]{@{}c@{}}Demand\\ assignment\end{tabular}} &  & \multirow{2}{*}{\begin{tabular}[c]{@{}c@{}}Multi-\\ period\end{tabular}} &  & \multirow{2}{*}{\begin{tabular}[c]{@{}c@{}}Freight\\ demand\end{tabular}} &  & \multirow{2}{*}{\begin{tabular}[c]{@{}c@{}}Network\\ application\end{tabular}} &  & \multirow{2}{*}{Solution method} \\
                            &  &                                                                                   &  &                                                                              &  &                                                                          &  &                                                                           &  &                                                                                &  &                                  \\ \hline
\cite{L18}                         &  & service                                                                         &  &                                                                              &  &                                                                          &  &                                                                           &  & \checkmark                                                                              &  & Commercial solver                \\
\cite{L5}                         &  & cost                                                                              &  &                                                                              &  &                                                                          &  &                                                                           &  & \checkmark                                                                              &  & Branch-and-bound                 \\
\cite{L3}                        &  & cost                                                                              &  &                                                                              &  &                                                                          &  &                                                                           &  & \checkmark                                                                              &  & Branch-and-cut                   \\
\cite{L4}                        &  & cost \& service                                                                 &  & \checkmark                                                                            &  &                                                                          &  &                                                                           &  &\checkmark                                                                              &  & Column generation                \\
\cite{L19}                         &  & profit                                                                            &  & \checkmark                                                                            &  &                                                                          &  &                                                                           &  &                                                                                &  & Branch-and-cut                   \\
\cite{L9}                          &  & profit                                                                            &  & \checkmark                                                                            &  &                                                                          &  &                                                                           &  & \checkmark                                                                              &  & Branch-and-cut                   \\
\cite{L1}                         &  & service                                                                         &  & \checkmark                                                                            &  &                                                                          &  &                                                                           &  & \checkmark                                                                             &  & Particle swarm                   \\
\cite{L16}                        &  & cost                                                                              &  &                                                                              &  & \checkmark                                                                        &  &                                                                           &  & \checkmark                                                                              &  & Commercial solver                \\
\cite{L20}                         &  & cost \& service                                                                             &  & \checkmark                                                                            &  &                                                                          &  &                                                                           &  &                                                                                &  & Lagrangian relaxation            \\
\cite{L12}                         &  & cost \& service                                                                 &  & \checkmark                                                                            &  &                                                                          &  &                                                                           &  & \checkmark                                                                             &  & Simulated annealing              \\
\cite{L8}                         &  & cost \& service                                                                 &  & \checkmark                                                                            &  &                                                                          &  &                                                                           &  & \checkmark                                                                             &  & Commercial solver              \\
\cite{L15}                        &  & profit                                                                            &  & \checkmark                                                                            &  &                                                                          &  &                                                                           &  & \checkmark                                                                              &  & Branch-and-bound                 \\
this paper                  &  & profit                                                                            &  & \checkmark                                                                            &  & \checkmark                                                                        &  & \checkmark                                                                        &  & \checkmark                                                                              &  & Column generation                \\ \hline
\end{tabular}
}
\label{review}
\end{table}

The line planning literature generally considers demand as static within a given planning horizon, with only a few recent papers considering the real-life dynamics of transport systems where demand can vary significantly throughout the day. For instance, \cite{L19} offer a schedule-based approach for intercity bus line planning, integrating dynamic demand and network planning. \cite{L16} present a multi-period line planning problem for urban transportation systems, addressing resource constraints and dynamic transfers. \cite{L20} propose a two-stage robust optimization model for railway line planning, considering passenger demand uncertainty. Additionally, in \cite{L12}, a bi-level programming model based on Stackelberg game theory optimizes line planning and flow distribution in high-speed railway networks, considering fluctuations in demand over time. 

Despite advancements in understanding passenger demand dynamics, the integration of freight demand into line planning has received relatively little attention. In scenarios where passenger demand varies, integrating freight transportation to allocate system capacity flexibly presents additional opportunities and complexities. This combined approach can potentially enhance the efficiency of infrastructure and resource utilization, thereby optimizing the overall performance of the railway network.

\subsection{Our contribution}
We have summarized the key features of the literature closely related to our study in Table \ref{review}. Our primary motivation is to propose a multi-period line planning approach that integrates passenger and freight transportation. To our knowledge, our work is the first attempt in this direction. Our approach is innovative in addressing time-varying passenger demand and incorporating the flexibility of freight demand at the strategic level of line planning. We employ column generation to dynamically generate paths, whereas most literature relies on a pre-defined set of paths. By leveraging advanced optimization techniques and heuristic algorithms, we can effectively solve large-scale instances in practical networks.

\section{Model formulation}
The classic line planning problem typically involves selecting a set of lines from a predefined line pool and determining their frequencies to efficiently meet given demand while respecting capacity constraints and operational objectives. This paper extends this problem to an integrated passenger-freight railway system, which introduces several additional considerations.

To take the different features of passenger demand and freight demand into account, our line planning model considers multiple periods, which represent, for example, peak and off-peak hours on a given day. Passenger demand is specified on a per-period level, while freight demand is specified for the entire day and can be served in all periods. Passengers and freight can only be assigned to paths with travel time below predefined thresholds. Not all demand needs to be served, but serving demand generates revenue for the operator. Each allowed line in the line pool specifies the type of rolling stock, or mode, used to operate the line, for example, a dedicated passenger train or a train that has capacity for both passengers and freight. Train capacity constraints are enforced at a disaggregated per-link level, necessitating detailed tracing of paths taken by passengers and freight. Finally, the objective is to maximize the total profit of railway operators, calculated as revenue from serving passengers and freight minus operational costs of the lines.

In the following sections, we discuss all components of the proposed line planning problem in detail,  including the underlying assumptions made and how we model it as an MIP.

\subsection{Railway network, periods and demand}
The railway network is represented as a graph $G=(V,E)$, where $V$ denotes the stations indexed by $i,j$, and the set $E$ defines the tracks as direct connections between stations. A subset of $V$ known as terminals, where lines can either start or end, is also present. 

Our line planning problem spans multiple discrete periods within a planning horizon, collectively forming the set $\mathcal{T} = \{1,2,...,T\}$ indexed by $t$. Lines can be installed to serve a set of passengers demand $P$ and freight demand $F$. Every $p\in P$ is defined by an origin $o_p\in V$, a destination $d_p\in V$, a period $t_p\in \mathcal{T}$ indicating in which period the passenger wishes to be served, a quantity $q_p$ and a unit revenue $\varphi_p$. Similarly, every $f\in F$ is defined by an origin $o_f\in V$, a destination $d_f\in V$, a quantity $q_f$ and a unit revenue $\varphi_f$. There is no period associated with freight demand since it is assumed freight can be served anytime during the planning horizon. 

\subsection{Line pool}
The concept of a line is fundamental in modeling the line planning problem. Trains on the railway network follow a predetermined sequence of non-repeating stops via tracks, forming a service path called a line. The set $\mathcal{K}$ contains all possible service routes indexed by $k$. Each line is also associated with a period $t\in \mathcal{T}$, indicating the period during which the line is operated, and a mode $m \in\mathcal{M}$, indicating the configuration and type of carriages of the line, which determines the capacity available for passengers and freight. A line is represented as a quadruple $(k,t,m)$,  and all potential lines constitute the line pool, denoted as $\mathcal{L}$. The operational cost of $l \in \mathcal{L}$ is denoted as $\xi_l$. For convenience, $\mathcal{L}(i)$ indicates the set of lines whose routes include station $i \in V$, $\mathcal{L}(e)$ indicates the set of lines that use track $e \in E$, and $\mathcal{L}(t)$ indicates the set of lines which operated in period $t \in \mathcal{T}$. 

\subsection{Passenger and freight paths}
To accurately trace the paths taken by passengers and freight, we use the change\&go-network, which was previously defined in \cite{L24,L25,L23}. We first describe the standard change\&go-network and then explain how it can be extended to account for different periods and modes. 

\subsubsection*{Standard change\&go-network}
The change\&go-network, denoted as $\mathcal{G=(V,A)}$, is constructed based on a physical graph $G$ as well as the line pool $\mathcal{L}$. Generally, the node set $\mathcal{V}$ consists of a set of travel nodes representing station-line pairs, denoted as  $\mathcal{V}_t = \{ (i,l) \in V \times \mathcal{L} : l \in \mathcal{L}(i)$\}, and station nodes introducing for each station $v \in V$, denoted as  $\mathcal{V}_s$. The arc set $\mathcal{A}$ consists of travel arcs that connect travel nodes along each line route, forming set $\mathcal{A}_t$, and transfer arcs that connect each station node to all corresponding travel nodes, forming set $\mathcal{A}_s$.

\begin{figure}[htbp]
\centering
\begin{subfigure}[b]{0.49\textwidth}
\centering
\begin{tikzpicture}[scale=1.5, every node/.style={scale=1.5}, >=stealth]

\coordinate (C) at (0,0);
\coordinate (D) at (2,0);
\coordinate (A) at (1,2.6);
\coordinate (B) at (1,1.3);


\draw (A) -- (B);
\draw (B) -- (C);
\draw (D) -- (B);
\draw[dashed] ([xshift=0 cm, yshift=0.5cm] A) -- (A);
\draw[dashed] (-0.305, -0.395) -- (C);
\draw[dashed]  (2.305, -0.395) -- (D);

\def\radius{0.3}
\draw[fill=white, draw=black] (A) circle (\radius) node {$a$};
\draw[fill=white, draw=black] (D) circle (\radius) node {$d$};
\draw[fill=white, draw=black] (B) circle (\radius) node {$b$};
\draw[fill=white, draw=black] (C) circle (\radius) node {$c$};

\draw[thick,->] ([xshift=-0.4cm, yshift=0cm] B) -- ([xshift=-0.1cm, yshift=0.4cm] C);
\draw[thick,] ([xshift=-0.4cm, yshift=1.25cm] B) -- ([xshift=-0.4cm, yshift=0cm] B);
\node [scale=0.7] at ([xshift=-1.3cm, yshift=0cm] B) {Line route 1};
\draw[thick,->] ([xshift=0.4cm, yshift=0cm] B) -- ([xshift=0.1cm, yshift=0.4cm] D);
\draw[thick] ([xshift=0.4cm, yshift=1.25cm] B) -- ([xshift=0.4cm, yshift=0cm] B) ;
\node [scale=0.7] at ([xshift=1.3cm, yshift=0cm] B) {Line route 2};
\end{tikzpicture}

\caption{Small railway network}
\label{cga}
\end{subfigure}
\begin{subfigure}[b]{0.49\textwidth}
\centering
\begin{tikzpicture}[scale=1.3, every node/.style={scale=1.3}, >=stealth]
\coordinate (C1) at (0,0);
\coordinate (D2) at (3,0);
\coordinate (A1) at (0,2.6);
\coordinate (A2) at (3,2.6);
\coordinate (B1) at (0,1.3);
\coordinate (B2) at (3,1.3);
\coordinate (TB) at (1.5,1.3);
\coordinate (TA) at (1.5,2.6);

\draw[thick,->] (A1) -- ([xshift=0cm, yshift=0.2cm] B1);
\draw[thick] (B1) -- ([xshift=0cm, yshift=0.2cm] C1);
\draw[thick,->] (A2) -- ([xshift=0cm, yshift=0.2cm] B2);
\draw[thick,->] (B2) -- ([xshift=0cm, yshift=0.2cm] D2);
\draw[dashed] ([xshift=0cm, yshift=0.5cm] A1) -- (A1);
\draw[dashed] ([xshift=0cm, yshift=0.5cm] A2) -- (A2);
\draw[dashed] (C1) -- ([xshift=0cm, yshift=-0.5cm]  C1);
\draw[dashed] (D2) -- ([xshift=0cm, yshift=-0.5cm]  D2);

\draw[->, thick, black!30] ([xshift=0.2cm, yshift=0.15cm]A1) .. controls (0.4,3) and (0.95,3) .. ([xshift=-0.2cm, yshift=0.15cm] TA);
\draw[->, thick, black!30] ([xshift=-0.2cm, yshift=-0.15cm] TA) .. controls (0.95,2.2) and (0.4,2.2) .. ([xshift=0.2cm, yshift=-0.15cm] A1);
\draw[->, thick, black!30] ([xshift=0.2cm, yshift=0.15cm]B1) .. controls (0.4,1.7) and (0.95,1.7) .. ([xshift=-0.2cm, yshift=0.15cm] TB);
\draw[->, thick, black!30] ([xshift=-0.2cm, yshift=-0.15cm] TB) .. controls (0.95,0.9) and (0.4,0.9) .. ([xshift=0.2cm, yshift=-0.15cm] B1);

\draw[->, thick, black!30] ([xshift=1.7cm, yshift=0.15cm]A1) .. controls (1.9,3) and (2.45,3) .. ([xshift=1.3cm, yshift=0.15cm] TA);
\draw[->, thick, black!30] ([xshift=1.3cm, yshift=-0.15cm] TA) .. controls (2.45,2.2) and (1.9,2.2) .. ([xshift=1.7cm, yshift=-0.15cm] A1);
\draw[->, thick, black!30] ([xshift=1.7cm, yshift=0.15cm]B1) .. controls (1.9,1.7) and (2.45,1.7) .. ([xshift=1.3cm, yshift=0.15cm] TB);
\draw[->, thick, black!30] ([xshift=1.3cm, yshift=-0.15cm] TB) .. controls (2.45,0.9) and (1.9,0.9) .. ([xshift=1.7cm, yshift=-0.15cm] B1);

\def\radius{0.2}
\draw[fill=black, draw=black] (C1) circle (\radius) node[xshift=-0.8cm] {$(c,l_1)$};
\draw[fill=black, draw=black] (D2) circle (\radius) node[xshift=0.8cm] {$(d,l_2)$};
\draw[fill=black, draw=black] (A1) circle (\radius) node[xshift=-0.8cm] {$(a,l_1)$};
\draw[fill=black, draw=black] (A2) circle (\radius) node[xshift=0.8cm] {$(a,l_2)$};
\draw[fill=black, draw=black] (B1) circle (\radius) node[xshift=-0.8cm] {$(b,l_1)$};
\draw[fill=black, draw=black] (B2) circle (\radius) node[xshift=0.8cm] {$(b,l_2)$};
\draw[fill=white, draw=black] (TA) circle (\radius) node[yshift=0.4cm] {$a$};
\draw[fill=white, draw=black] (TB) circle (\radius) node[yshift=0.4cm]  {$b$};

\node[draw, circle, fill=black, inner sep=1pt] at (0,-0.8) {};
\node [scale=0.5] at (1,-0.8) {Travel Node};
\node[draw, circle, inner sep=1pt] at (2,-0.8) {};
\node [scale=0.5] at (2.9,-0.8) {Station Node};

\draw[->] (-0.50,-1.1) -- (0.3,-1.1);
\node [scale=0.5] at (1.1,-1.1) {Travel Arc};
\draw[->,thick, black!30] (1.9,-1.1) -- (2.6,-1.1);
\node [scale=0.5] at (3.3,-1.1) {Transfer Arc};

\end{tikzpicture}
\caption{Standard change\&go-network}
\label{cgb}
\end{subfigure}
\caption{Example of a standard change\&go-network}
\label{cg1}
\end{figure}
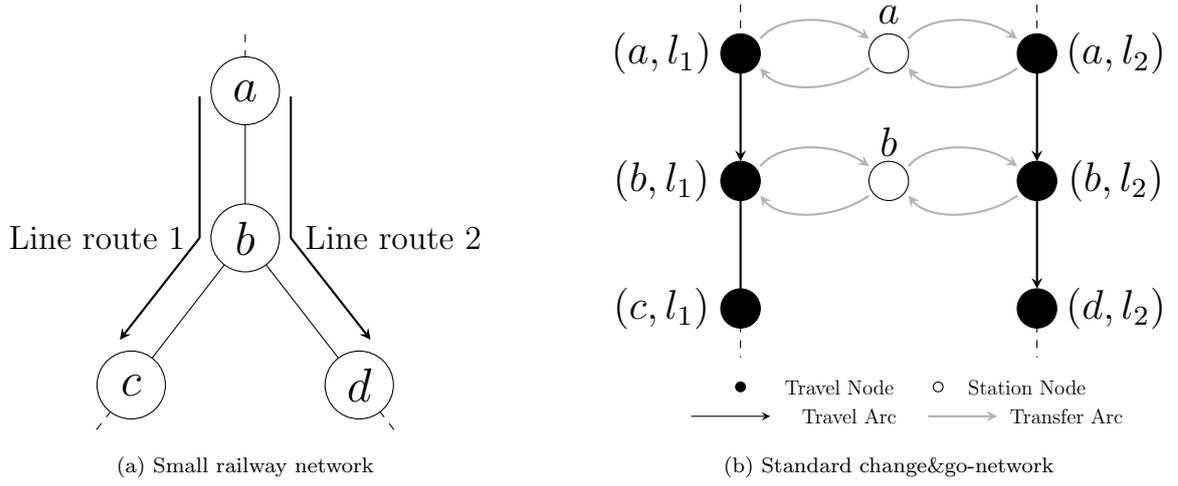

To illustrate, consider the small railway network with four stations and three tracks shown in Figure \ref{cga}, accompanied by the line pool $\mathcal{L} = \{l_1, l_2\} = \{(a,b,c), (a,b,d)\}$ disregarding specific periods and modes to illustrate a standard change\&go-network. As shown in Figure \ref{cgb}, we get the standard change\&go-network with the travel node set $\mathcal{V}_t = \{ (a, l_1),(b, l_1),(c, l_1),(a, l_2),(b, l_2),$ $(d, l_2)\}$, where each element is interconnected via travel arcs based on the line routes, and station node set $\mathcal{V}_s = \{ a, b\}$, connected to travel nodes within the same station via transfer arcs. 

\subsubsection*{Period-extended change\&go-network}

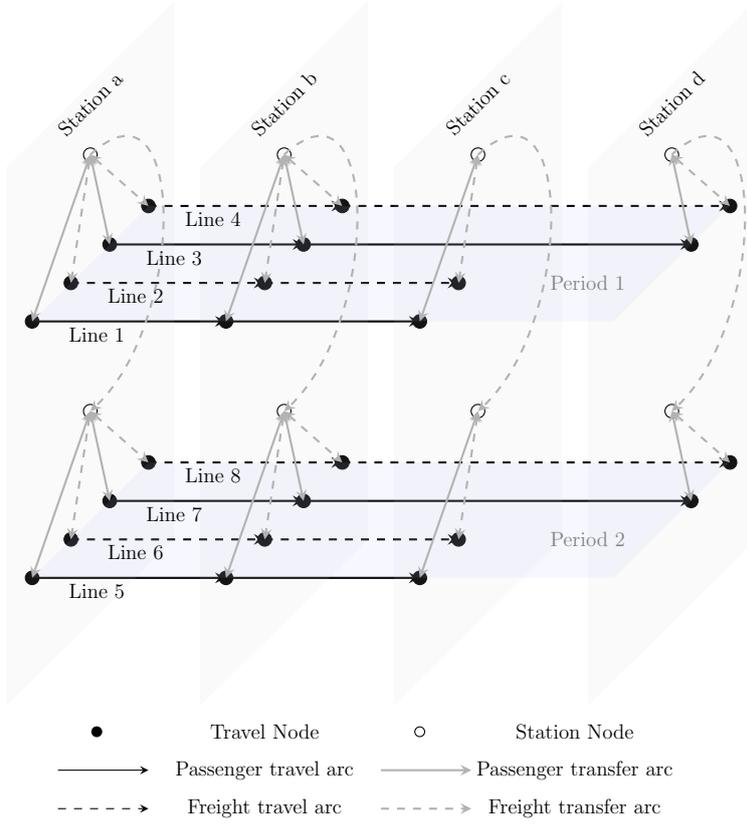
\begin{figure}[htbp]
\centering
\begin{tikzpicture}[scale=1.7, every node/.style={scale=1.3}, >=stealth]

\coordinate (11) at (0,0);
\coordinate (12) at (1.5,0);
\coordinate (13) at (3,0);
\coordinate (14) at (4.5,0);

\coordinate (21) at (0.3,0.3);
\coordinate (22) at (1.8,0.3);
\coordinate (23) at (3.3,0.3);
\coordinate (24) at (4.8,0.3);

\coordinate (31) at (0.6,0.6);
\coordinate (32) at (2.1,0.6);
\coordinate (33) at (3.6,0.6);
\coordinate (34) at (5.1,0.6);

\coordinate (41) at (0.9,0.9);
\coordinate (42) at (2.4,0.9);
\coordinate (43) at (3.9,0.9);
\coordinate (44) at (5.4,0.9);

\coordinate (51) at (0,2);
\coordinate (52) at (1.5,2);
\coordinate (53) at (3,2);
\coordinate (54) at (4.5,2);

\coordinate (61) at (0.3,2.3);
\coordinate (62) at (1.8,2.3);
\coordinate (63) at (3.3,2.3);
\coordinate (64) at (4.8,2.3);

\coordinate (71) at (0.6,2.6);
\coordinate (72) at (2.1,2.6);
\coordinate (73) at (3.6,2.6);
\coordinate (74) at (5.1,2.6);

\coordinate (81) at (0.9,2.9);
\coordinate (82) at (2.4,2.9);
\coordinate (83) at (3.9,2.9);
\coordinate (84) at (5.4,2.9);

\foreach \i in {1,2,5,6} {
    \draw[black, fill=black] (\i 1) circle (1.5pt);
    \draw[black, fill=black] (\i 2) circle (1.5pt);
    \draw[black, fill=black] (\i 3) circle (1.5pt);
}

\foreach \i in {3,4,7,8} {
    \draw[black, fill=black] (\i 1) circle (1.5pt);
    \draw[black, fill=black] (\i 2) circle (1.5pt);
    \draw[black, fill=black] (\i 4) circle (1.5pt);
}

\foreach \i in {1,5} {
    \draw[thick,->, black] (\i 1) -- (\i 2);
    \draw[thick,->, black] (\i 2) -- (\i 3);
}

\foreach \i in {2,6} {
    \draw[dashed, thick,->, black] (\i 1) -- (\i 2);
    \draw[dashed, thick,->, black] (\i 2) -- (\i 3);
}

\foreach \i in {3,7} {
    \draw[thick,->, black] (\i 1) -- (\i 2);
    \draw[thick,->, black] (\i 2) -- (\i 4);
}

\foreach \i in {4,8} {
    \draw[dashed, thick,->, black] (\i 1) -- (\i 2);
    \draw[dashed, thick,->, black] (\i 2) -- (\i 4);
}

\fill[blue!30, opacity=0.1] (11) -- (14) -- (44) -- (41) -- cycle;
\fill[blue!30, opacity=0.1] (51) -- (54) -- (84) -- (81) -- cycle;

\fill[gray!40, opacity=0.1] ([xshift=-0.2cm, yshift=-1cm] 11) -- ([xshift=0.2cm, yshift=-0.6cm] 41) -- ([xshift=0.2cm, yshift=1.6cm] 81) -- ([xshift=-0.2cm, yshift=1.2cm] 51) -- cycle;
\fill[gray!40, opacity=0.1] ([xshift=-0.2cm, yshift=-1cm] 12) -- ([xshift=0.2cm, yshift=-0.6cm] 42) -- ([xshift=0.2cm, yshift=1.6cm] 82) -- ([xshift=-0.2cm, yshift=1.2cm] 52) -- cycle;
\fill[gray!40, opacity=0.1] ([xshift=-0.2cm, yshift=-1cm] 13) -- ([xshift=0.2cm, yshift=-0.6cm] 43) -- ([xshift=0.2cm, yshift=1.6cm] 83) -- ([xshift=-0.2cm, yshift=1.2cm] 53) -- cycle;
\fill[gray!40, opacity=0.1] ([xshift=-0.2cm, yshift=-1cm] 14) -- ([xshift=0.2cm, yshift=-0.6cm] 44) -- ([xshift=0.2cm, yshift=1.6cm] 84) -- ([xshift=-0.2cm, yshift=1.2cm] 54) -- cycle;

\foreach \i in {1,2,3,4} {
    \draw[black] ([xshift=0.15cm, yshift=1cm] 6\i) circle (1.5pt);
}
\draw[<->, thick, black!30] ([xshift=0.15cm, yshift=1cm] 61) -- (51);
\draw[<->, thick, black!30] ([xshift=0.15cm, yshift=1cm] 61) -- (71);
\draw[<->, thick, black!30] ([xshift=0.15cm, yshift=1cm] 62) -- (52);
\draw[<->, thick, black!30] ([xshift=0.15cm, yshift=1cm] 62) -- (72);
\draw[<->, thick, black!30] ([xshift=0.15cm, yshift=1cm] 63) -- (53);
\draw[<->, thick, black!30] ([xshift=0.15cm, yshift=1cm] 64) -- (74);

\draw[dashed, <->, thick, black!30] ([xshift=0.15cm, yshift=1cm] 61) -- (61);
\draw[dashed,<->, thick, black!30] ([xshift=0.15cm, yshift=1cm] 61) -- (81);
\draw[dashed,<->, thick, black!30] ([xshift=0.15cm, yshift=1cm] 62) -- (62);
\draw[dashed,<->, thick, black!30] ([xshift=0.15cm, yshift=1cm] 62) -- (82);
\draw[dashed,<->, thick, black!30] ([xshift=0.15cm, yshift=1cm] 63) -- (63);
\draw[dashed,<->, thick, black!30] ([xshift=0.15cm, yshift=1cm] 64) -- (84);

\foreach \i in {1,2,3,4} {
    \draw[black] ([xshift=0.15cm, yshift=1cm] 2\i) circle (1.5pt);
}
\draw[<->, thick, black!30] ([xshift=0.15cm, yshift=1cm] 21) -- (11);
\draw[<->, thick, black!30] ([xshift=0.15cm, yshift=1cm] 21) -- (31);
\draw[<->, thick, black!30] ([xshift=0.15cm, yshift=1cm] 22) -- (12);
\draw[<->, thick, black!30] ([xshift=0.15cm, yshift=1cm] 22) -- (32);
\draw[<->, thick, black!30] ([xshift=0.15cm, yshift=1cm] 23) -- (13);
\draw[<->, thick, black!30] ([xshift=0.15cm, yshift=1cm] 24) -- (34);
\draw[dashed,<->, thick, black!30] ([xshift=0.15cm, yshift=1cm] 21) -- (21);
\draw[dashed,<->, thick, black!30] ([xshift=0.15cm, yshift=1cm] 21) -- (41);
\draw[dashed,<->, thick, black!30] ([xshift=0.15cm, yshift=1cm] 22) -- (22);
\draw[dashed,<->, thick,black!30] ([xshift=0.15cm, yshift=1cm] 22) -- (42);
\draw[dashed,<->, thick,black!30] ([xshift=0.15cm, yshift=1cm] 23) -- (23);
\draw[dashed,<->, thick, black!30] ([xshift=0.15cm, yshift=1cm] 24) -- (44);

\draw[dashed, ->, thick, black!30] ([xshift=0.15cm, yshift=1cm] 61) .. controls (1.2,4) and (1.2,2) .. ([xshift=0.15cm, yshift=1cm] 21);
\draw[dashed, ->, thick, black!30] ([xshift=0.15cm, yshift=1cm] 62) .. controls (2.7,4) and (2.7,2) .. ([xshift=0.15cm, yshift=1cm] 22);
\draw[dashed, ->, thick, black!30] ([xshift=0.15cm, yshift=1cm] 63) .. controls (4.2,4) and (4.2,2) .. ([xshift=0.15cm, yshift=1cm] 23);
\draw[dashed, ->, thick, black!30] ([xshift=0.15cm, yshift=1cm] 64) .. controls (5.7,4) and (5.7,2) .. ([xshift=0.15cm, yshift=1cm] 24);

\node[gray,scale=0.5] at (4.3,2.3) {Period 1};
\node[gray,scale=0.5] at (4.3,0.3) {Period 2};
\node[scale=0.5] at (0.5,-0.1) {Line 5};
\node[scale=0.5] at (0.8,0.2) {Line 6};
\node[scale=0.5] at (1.1,0.5) {Line 7};
\node[scale=0.5] at (1.4,0.8) {Line 8};
\node[scale=0.5] at (0.5,1.9) {Line 1};
\node[scale=0.5] at (0.8,2.2) {Line 2};
\node[scale=0.5] at (1.1,2.5) {Line 3};
\node[scale=0.5] at (1.4,2.8) {Line 4};

\node[scale=0.5, rotate=45] at (0.45,3.7) {Station a};
\node[scale=0.5, rotate=45] at (1.95,3.7) {Station b};
\node[scale=0.5, rotate=45] at (3.45,3.7) {Station c};
\node[scale=0.5, rotate=45] at (4.95,3.7) {Station d};

\node[draw, circle, fill=black, inner sep=1pt] at (0.5,-1.2) {};
\node [scale=0.5] at (1.8,-1.2) {Travel Node};
\node[draw, circle, inner sep=1pt] at (3,-1.2) {};
\node [scale=0.5] at (4.2,-1.2) {Station Node};

\draw[->] (0.2,-1.5) -- (0.9,-1.5);
\node [scale=0.5] at (1.8,-1.5) {Passenger travel arc};
\draw[->,thick, black!30] (2.7,-1.5) -- (3.4,-1.5);
\node [scale=0.5] at (4.2,-1.5) {Passenger transfer arc};

\draw[dashed, ->] (0.2,-1.8) -- (0.9,-1.8);
\node [scale=0.5] at (1.8,-1.8) {Freight travel arc};
\draw[dashed, ->,thick, black!30] (2.7,-1.8) -- (3.4,-1.8);
\node [scale=0.5] at (4.2,-1.8) {Freight transfer arc};

\end{tikzpicture}
\caption{An example on period-extended change\&go-network}
\label{cg2}
\end{figure}

Building upon the standard change\&go-network, we enhance it to accommodate the multi-period line planning problem, referred to as a period-extended change\&go-network. Based on the network in Figure \ref{cga}, we now introduce two periods and two modes: dedicated passenger trains and dedicated freight trains. The line pool expands to include eight lines, detailed in Table \ref{lp}. The corresponding period-extended change\&go-network is shown in Figure \ref{cg2}.

\begin{table}[htbp]
\renewcommand{\arraystretch}{0.8}
\centering
\caption{The line pool of the small network considering two periods and two modes}
\adjustbox{max width=1\textwidth}{
\begin{tabular}{clclclc}
\hline
Period                   &  & Line   &  & Route &  & Mode            \\ \hline
\multirow{4}{*}{Period 1} &  & Line 1 &  & (a,b,c)       &  & Passenger train \\
                         &  & Line 2 &  & (a,b,c)       &  & Freight train   \\
                         &  & Line 3 &  & (a,b,d)       &  & Passenger train \\
                         &  & Line 4 &  & (a,b,d)       &  & Freight train   \\ \hline
\multirow{4}{*}{Period 2} &  & Line 5 &  & (a,b,c)       &  & Passenger train \\
                         &  & Line 6 &  & (a,b,c)       &  & Freight train   \\
                         &  & Line 7 &  & (a,b,d)       &  & Passenger train \\
                         &  & Line 8 &  & (a,b,d)       &  & Freight train   \\ \hline
\end{tabular}
}
\label{lp}
\end{table}

Figure \ref{cg2} shows that travel nodes are generated in two levels according to their belonging period. Notably, station nodes are now introduced for every station within each period, denoted as $\mathcal{V}_s = \{(i,t) \in V \times \mathcal{T} \}$. Travel arcs connect travel nodes, while transfer arcs link station nodes to travel nodes within the same station and period. Moreover, there are transfer arcs between two consecutive periods, but only for freight. Therefore, freight paths can comprise multiple periods, while passenger paths are enforced to be contained in a single period. 

 The travel time of a travel arc $\{(i,l),(j,l)\}$ is calculated as the sum of the running time from station $i$ to $j$ and the stopping time. Transfer times are accounted for in transfer arcs from a station node to a travel node and across periods, while the travel time from a travel node to a station node is assumed to be zero.

\subsubsection*{Paths}
Given passenger demand $p \in P$, the potential paths are all paths from travel nodes of station $o_p$ in period $t_p$ to travel nodes of station $d_p$ in period $t_p$. For freight demand $f \in F$, the potential paths are all paths from travel nodes of station $o_f$ to travel nodes of station $d_f$. To ensure service quality, we only consider simple paths whose duration does not exceed a designated threshold. The threshold, denoted as  $\tilde{\tau_p}$ for $p \in P$ and $\tilde{\tau_f}$ for  $f \in F$, is determined by the minimal travel time on change\&go-network and a multiplication factor for each demand. Freight demand typically has a larger factor than passengers, as freight can tolerate longer detours. Additionally, for the same type of demand, the greater the minimal travel time, the longer the acceptable detour.

For convenience, the following notations are introduced. Each travel arc $a \in \mathcal{A}_t$ is associated with a line $l_a\in \mathcal{L}$ and characterized by a passenger capacity $c_a$ and a freight capacity $c_a^{'}$, determined by the mode of the line $l_a$. The path set for passenger demand $p \in P$ is denoted as $\mathcal{R}_p$, and $\mathcal{R}_p(a)$ represents the subset using arc $a \in \mathcal{A}$. The path pool for freight demand $f \in F$ is denoted as set $\mathcal{R}_f$, and $\mathcal{R}_f(a)$ denotes the paths using arc $a \in \mathcal{A}$.

\subsection{Formulation}
The multi-period line planning problem can be modeled as an MIP model using three kinds of variables:
\begin{enumerate}[i)]
    \item non-negative integer variables $x_l$ to indicate the frequency of line $l \in \mathcal{L}$, 
    \item non-negative continuous variables $y_r$ to indicate the amount of passengers traveling using route $r \in \mathcal{R}_p, p \in P$, and
    \item non-negative continuous variables $z_r$ for the amount of freight on route $r \in \mathcal{R}_f, f \in F$.
\end{enumerate}

\begin{maxi!}
    {x_l,y_r,z_r} 
    { \sum_{p \in P}  \varphi_p \sum_{r\in \mathcal{R}_p} y_r  + \sum_{f \in F}   \varphi_f \sum_{r\in \mathcal{R}_f}  z_r - \sum_{l \in \mathcal{L}} \xi_l x_l} 
    {}{} \label{obj}
    \addConstraint{\sum_{r\in \mathcal{R}_p} y_r}{\le q_p, \quad \forall p \in P} \label{c1}
    \addConstraint{\sum_{r\in \mathcal{R}_f} z_r}{\le q_f, \quad \forall f \in F} \label{c2}
    \addConstraint{\sum_{r\in \mathcal{R}_p(a)}\sum_{p \in P} y_r - c_a x_{l_a}}{\le 0, \quad \forall a\in \mathcal{A}_t} \label{c3}
    \addConstraint{\sum_{r\in \mathcal{R}_f(a)}\sum_{f \in F} z_r - c_a^{'} x_{l_a}}{\le 0, \quad \forall a\in \mathcal{A}_t} \label{c4}
    \addConstraint{\sum_{l \in \mathcal{L}(e) \cap \mathcal{L}(t)} x_l }{\le \Lambda_e^t, \quad \forall e \in E, t \in \mathcal{T}} \label{c5}
    \addConstraint{x_l}{\in \mathbb{N}, \quad \forall l\in \mathcal{L}}\label{c6}
    \addConstraint{y_r}{\ge 0, \quad \forall p \in P,  r\in \mathcal{R}_p} \label{c7}
    \addConstraint{z_r}{\ge 0, \quad \forall  f \in F, r\in \mathcal{R}_f.} \label{c8}
\end{maxi!}

The objective of the model \ref{obj} is to maximize profit, specified as the difference between revenue generated from ticket sales and line costs. Constraints \ref{c1} ensure that the number of routed passengers of demand $p \in P$ is at most $q_p$, and constraints \ref{c3} link the route variables with the line variables and guarantee sufficient capacity for the routed passengers. Constraints \ref{c2} and \ref{c4} serve the same purpose for the freight demand.  The maximal throughput capacity of track $e \in E$ in period $t \in \mathcal{T}$ is set to $\Lambda_e^t$ in constraints \ref{c5}.

\section{Solution approach}
For larger instances, the number of paths becomes enormous, making it impossible to solve the MIP model directly. Therefore, we opt for a column generation approach to solve the model's linear relaxation. Subsequently, we use heuristics to obtain an integer solution.

\subsection{Column generation}
To manage the vast number of possible paths, column generation is applied to generate paths for passengers and freight dynamically. We relax the MIP model to an LP, called the master problem (MP). An initial path pool is constructed by generating paths with minimal travel time for each demand. We then solve the restricted master problem (RMP) with a subset of paths. Using the dual variables obtained from the RMP, we solve the pricing problems to identify paths with positive reduced costs (the MP is a maximization problem). This iterative process continues until no further paths with positive reduced costs can be found. Next, we will explain the specific details of the implementation. 

\subsubsection{Pricing problems}
The dual variables from solving the RMP are denoted as $\alpha_p$ for passenger flow constraints (1b), $\beta_f$ for freight flow constraints (1c), $\delta_a$ for passenger capacity constraints (1d), and $\varepsilon_a$ for freight capacity constraints (1e). The RMP is always feasible because a solution with zero profit can always be found by setting all frequencies to zero. Therefore, we can always determine the values of the dual variables by solving the RMP.

The reduced cost ${{\overline{\tau }}_r}$ for variable $y_r$ with $r\in \mathcal{R}_p, p \in P$ satisfies
\begin{equation}
{{\overline \tau}_r}= \varphi_p - \alpha _p - \sum\limits_{a\in r}{\delta _{a}}= (\varphi_p - \alpha _{p}) - \sum\limits_{a\in r}{\delta _{a}}
\end{equation}
Therefore, maximizing reduced cost reduces to finding a least-cost simple path in change\&go-network concerning the non-negative costs $\delta _{a}$ on the arcs, such that ${{\overline{\tau }}_{r}} > 0$ and the total duration is no more than $\tilde{\tau_p}$. If the least-cost path has a negative reduced cost, it can be concluded that there is no path with a positive reduced cost. 

Similarly, the reduced cost ${{\overline{\gamma }}_{r}}$ for variable $z_r$ with $r \in \mathcal{R}_f, f \in F$, is:
\begin{equation}
{{\overline{\gamma }}_{r}}= \varphi_f - \beta_f - \sum\limits_{a \in r}\varepsilon_a = (\varphi_f - \beta_f)  - \sum\limits_{a \in r} \varepsilon_a
\end{equation}

Hence, the pricing problem for the $z$-variables is the same as for the $y$-variables, but with the arc costs $\varepsilon_a $ and duration bound $\tilde{\tau_f}$.

\subsubsection{Dijkstra and label-correcting}
Because the arc costs are non-negative in both pricing problems, they can be solved as the shortest-path problems with duration constraints. This problem can be addressed using the classic label-correcting algorithm. Since this can be time-consuming, we enhance efficiency by integrating Dijkstra algorithm to solve the shortest-path problem without duration constraints, corresponding to a relaxation of the pricing problem. If the resulting path has positive reduced costs but violates the duration constraint, we apply label-correcting to find the least-cost path satisfying the duration constraint. 

The shortest path for each demand, respecting the dual costs on arcs, can be found easily using the Dijkstra algorithm in polynomial time. For example, for passenger demand $p \in P$ with a positive value of $\varphi_p - \alpha _p$, we find the shortest path in change\&go-network with respect to the costs $\delta _a$ on arcs. For the shortest path found via Dijkstra, no valid path can be added to the model if its reduced cost is negative. Otherwise, we check if the path's duration meets the duration constraint; if it does, it is directly added to the path pool. If not, we proceed with the label-correcting algorithm to search the least-cost path whose duration is below the threshold.

In label-correcting, each path from the origin node to a given node is associated with a label containing the total cost and travel time. A path from the current node can be extended to another node only if there is a direct outgoing arc and the duration threshold is not exceeded after the extension. When two paths converge at the same node, the one with the smaller cost and shorter duration dominates. Throughout the algorithm's execution, we only consider non-dominated paths and extendable nodes. For more detailed information about the algorithm, please refer to \cite{A2} and \cite{A1}.

Transportation can start and end via different lines for each demand, resulting in multiple potential start and end nodes. However, these start and end nodes are associated with specific station nodes. Therefore, for passenger demand $p \in P$, the pricing problem can be solved by searching for the shortest path from station node $(o_p, t_p)$ to station node $(d_p, t_p)$. Similarly, for freight demand $f \in F$, the pricing problem can be solved by searching for the shortest paths from station node $(o_f, 1)$ to station node $(d_f, t)$, where $t \in \mathcal{T}$, and selecting the path with maximally reduced cost among them. The cost and duration of a path are calculated starting from the first routed travel node. The Dijkstra and label-correcting algorithms are single-source shortest-path algorithms, and the pricing network is the same for all OD pairs. Hence, for all passenger or freight demand originating from the same station node, the pricing problem can be consolidated and solved by executing the algorithm once.

\subsection{Obtaining integer solutions}

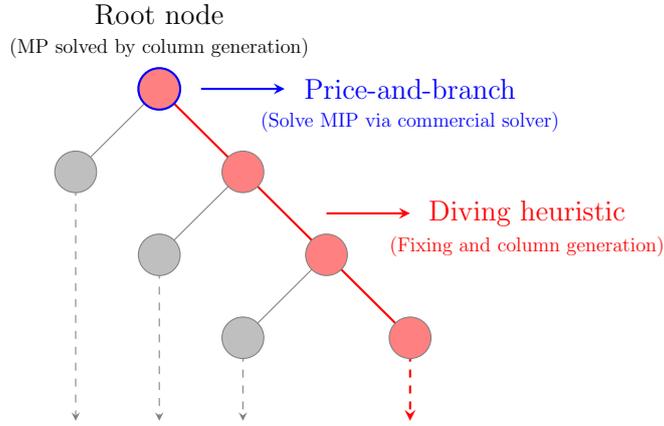
\begin{figure}[ht]
    \centering
    \begin{tikzpicture}[scale=1.1, every node/.style={scale=0.9}, >=stealth]

\coordinate (A) at (1,4);
\coordinate (B2) at (0,3);
\coordinate (B1) at (2,3);
\coordinate (C2) at (1,2);
\coordinate (C1) at (3,2);
\coordinate (D2) at (2,1);
\coordinate (D1) at (4,1);
\coordinate (E4) at (0,0);
\coordinate (E3) at (1,0);
\coordinate (E2) at (2,0);
\coordinate (E1) at (4,0);

\draw[red,thick] (A) -- (B1);
\draw[gray] (A) -- (B2);
\draw[red,thick] (B1) -- (C1);
\draw[gray] (B1) -- (C2);
\draw[red,thick] (C1) -- (D1);
\draw[gray] (C1) -- (D2);
\draw[red,thick,dashed,->] (D1) -- (E1);
\draw[gray,dashed,->]  (D2) -- (E2);
\draw[gray,dashed,->]  (C2) -- (E3);
\draw[gray,dashed,->]  (B2) -- (E4);

\def\radius{0.3}
\draw[fill=red!50, draw=blue, thick] (A) circle (0.25cm);
\draw[fill=red!50, draw=gray] (B1) circle (0.25cm);
\draw[fill=red!50, draw=gray] (C1) circle (0.25cm);
\draw[fill=red!50, draw=gray] (D1) circle (0.25cm);
\draw[fill=gray!50, draw=gray] (B2) circle (0.25cm);
\draw[fill=gray!50, draw=gray] (C2) circle (0.25cm);
\draw[fill=gray!50, draw=gray] (D2) circle (0.25cm);

\node [scale=1] at ([yshift=0.9cm] A) {Root node};
\node [scale=0.7] at ([yshift=0.5cm] A) {(MP solved by column generation)};

\draw[blue, thick, ->] ([xshift=0.5cm] A) -- ([xshift=1.5cm] A);
\node [blue, scale=1] at ([xshift=3cm] A) {Price-and-branch};
\node [blue, scale=0.7] at ([xshift=3cm, yshift=-0.4cm] A) {(Solve MIP via commercial solver)};

\draw[red, thick, ->] (3,2.5) -- (4,2.5);
\node [red, scale=1] at (5.4,2.5) {Diving heuristic};
\node [red, scale=0.7] at (5.4,2.1) {(Fixing and column generation)};
\end{tikzpicture}
    \caption{Illustration of two heuristics for obtaining an integer solution}
    \label{Fig alg}
\end{figure}

We now discuss two heuristics to find integer solutions for the proposed line planning problem. Initially, we determine an upper bound in both cases by solving the MP through column generation, as shown in Figure \ref{Fig alg}.

Our first heuristic, the restricted master heuristic or price-and-branch method, applies a commercial solver to the MIP model using all columns generated in solving the MP.

Our second heuristic involves a depth-first exploration of the branch-and-bound tree, commonly referred to as a diving heuristic in the literature. In this heuristic, two phases, namely fixing and column generation, are executed iteratively until an integer solution is obtained. After a column generation phase is terminated, the fractional path variable whose value is nearest to a non-zero integer will be rounded and fixed to that integer. Then, we resolve the RMP and start the column generation process again. We repeat these two steps until all variables in the solution are integer. 

To enhance the efficiency of the column generation phase in the diving algorithm, we terminate this phase once a pre-defined maximum number of iterations without improvement is reached. Secondly, during execution, we reduce the size of the change\&go-network by eliminating specific travel nodes and their adjacent arcs when the corresponding lines utilize tracks that have reached their maximum throughput capacity due to the fixed lines.

\section{Experimental design}
\begin{figure}[t]
    \centering    \includegraphics[width=.6\columnwidth]{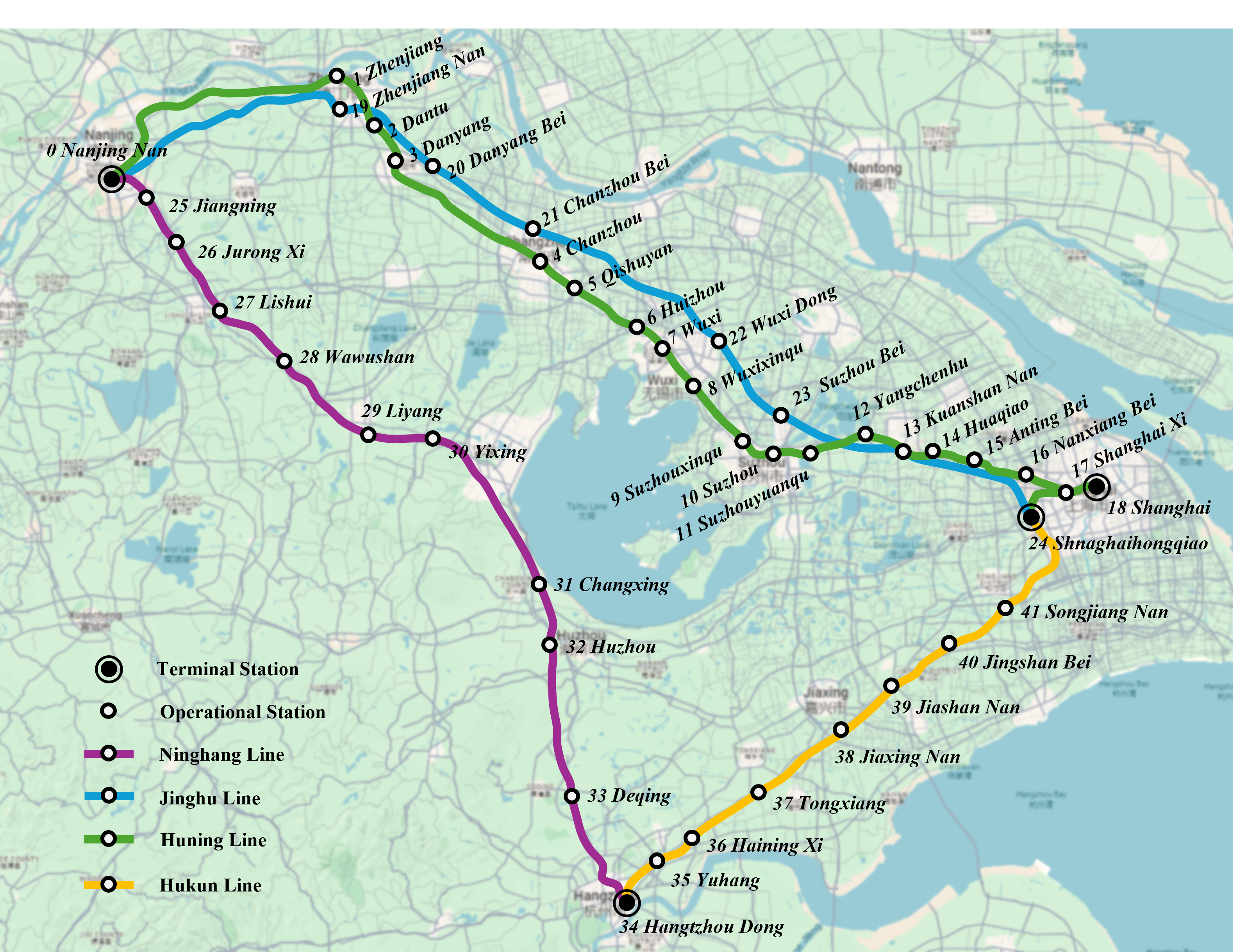}
    \caption{Large network from complete Chinese high-speed railway network}
    \label{full network}
\end{figure}

\begin{figure}[t]
    \centering    \includegraphics[width=.6\columnwidth]{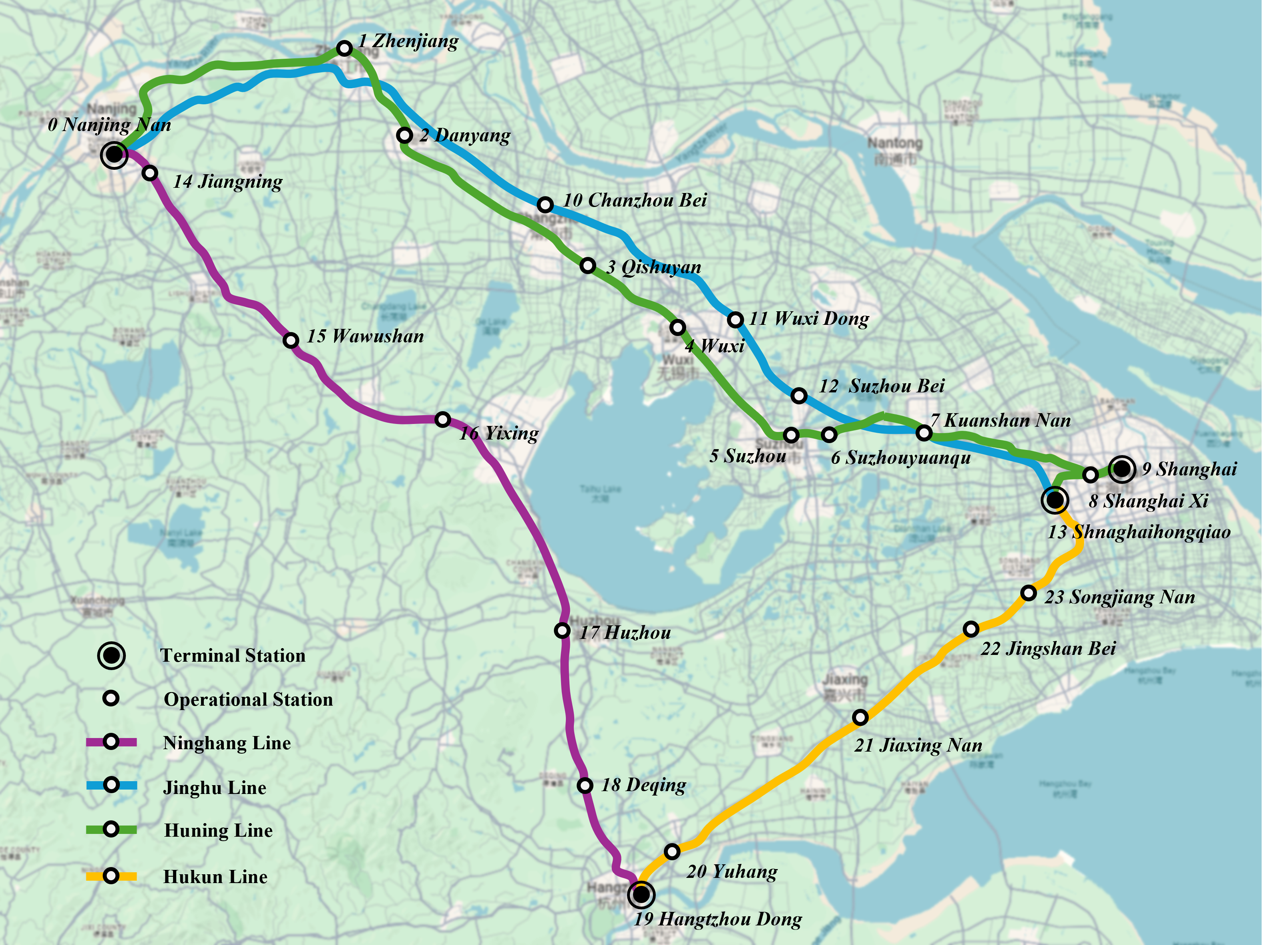}
    \caption{Medium network from reduced Chinese high-speed railway network}
    \label{reduced network}
\end{figure}

The experiments are conducted on two sets of instances: a larger set and a smaller set, both derived from a real-life Chinese high-speed railway network in the Yangtze River Delta region. The larger set of instances is generated from the large network, which includes four key high-speed railway lines (Ninghang, Jinghu, Huning, and Hukun) and comprises 42 stations, including four terminal stations, as shown in Figure~\ref{full network}. The lengths of the tracks are known. The smaller set of instances is created by selecting a subset of operational stations to form a medium network, as depicted in Figure \ref{reduced network}.

In all instances, the planning horizon spans from 8:00 to 20:00 and is divided into three distinct periods: Morning peak period (M-peak) from 8:00 to 12:00, off-peak hours from 12:00 to 16:00, and evening peak period (E-peak) from 16:00 to 20:00. 

A line can be performed using one of three modes: dedicated passenger train, dedicated freight train, or a mixed train with half of the capacity for carrying freight and the other half for passengers. Each instance considers three operational scenarios:
 \begin{enumerate}[i)]
     \item P: Serving only passengers, with the mode set including solely passenger train;
     \item P+F: Serving both passengers and freight, with the mode set including only passenger and freight trains;
     \item P+F+M: Serving both passengers and freight, with the mode set incorporating passenger, freight, and mixed trains. 
 \end{enumerate}

\subsection{Line pool generation}
For the line pool, we need to further specify the set of line routes. We pre-generate line routes including all simple paths between any two terminals. For the stopping pattern, we consider two stop schemes: stopping at every routed station and stopping only at routed terminals. However, in principle, any stopping pattern could be applied to customize the line pool to meet specific operational or logistical needs if necessary. Due to the identical network topology, both the medium and the large networks share the same set of line routes, totaling 124 routes each. 

\subsection{Demand generation}

We randomly generate passenger and freight demand for each OD pair during the planning horizon using pre-assumed expected values, as detailed in Appendix A. Assuming equal passenger demand during M-peak and E-peak periods, the generated passenger demand for each OD pair is allocated across different periods using the parameter $R$ to denote the ratio of passenger demand in the peak period versus the off-peak period. For example, $R=1$ corresponds to uniform demand, and $R=2$ implies that the demand in the peak periods is twice the demand in the off-peak period. In general, the proportions allocated to the M-peak, off-peak, and E-peak periods are $R/(2R+1)$, $1/(2R+1)$, and $R/(2R+1)$, respectively. For medium and large networks, we both generated 10 instances with different demand levels, as specified in Table \ref{Inst}. 

\begin{table}[t]
\renewcommand{\arraystretch}{0.8}
\centering
\caption{Description of medium and large instances}
\adjustbox{max width=0.9\textwidth}{
\begin{tabular}{ccccccccccccccccc}
\hline
\multirow{3}{*}{\centering \begin{tabular}[c]{@{}c@{}}Demand \\ level\end{tabular}} &  & \multicolumn{7}{c}{Medium instance}                             &  & \multicolumn{7}{c}{Large instance}                             \\ \cline{3-9} \cline{11-17} 
                          &  & \multirow{2}{*}{\begin{tabular}[c]{@{}c@{}}Number \\ of ODs\end{tabular}} &  & \multicolumn{5}{c}{Size of demand}   &  & \multirow{2}{*}{\begin{tabular}[c]{@{}c@{}}Number \\ of ODs\end{tabular}} &  & \multicolumn{5}{c}{Size of demand}   \\ \cline{5-9} \cline{13-17} 
                          &  &                          &  & Passenger &  & Freight &  & Total &  &                          &  & Passenger &  & Freight &  & Total \\ \cline{1-17}
1                         &  & 100                      &  & 300       &  & 100     &  & 400   &  & 180                      &  & 540       &  & 180     &  & 720   \\
2                         &  & 150                      &  & 450       &  & 150     &  & 600   &  & 360                      &  & 1080      &  & 360     &  & 1440  \\
3                         &  & 200                      &  & 600       &  & 200     &  & 800   &  & 540                      &  & 1620       &  & 540     &  & 2160   \\
4                         &  & 250                      &  & 750       &  & 250     &  & 1000  &  & 720                      &  & 2160      &  & 720     &  & 2880  \\
5                         &  & 300                      &  & 900      &  & 300     &  & 1200   &  & 900                      &  & 2700       &  & 900     &  & 3600  \\
6                         &  & 350                      &  & 1050      &  & 350     &  & 1400  &  & 1080                     &  & 3240      &  & 1080    &  & 4320  \\
7                         &  & 400                      &  & 1200       &  & 400     &  & 1600   &  & 1260                     &  & 3780       &  & 1260    &  & 5040  \\
8                         &  & 450                      &  & 1350      &  & 450     &  & 1800  &  & 1440                     &  & 4320      &  & 1440    &  & 5760  \\
9                         &  & 500                      &  & 1500       &  & 500     &  & 2000   &  & 1620                     &  & 4860       &  & 1620    &  & 6480  \\
10                        &  & 552                      &  & 1656      &  & 552     &  & 2208  &  & 1722                     &  & 5166      &  & 1722    &  & 6888  \\ \hline
\end{tabular}
}
\label{Inst}
\end{table}

\subsection{Remaining parameter settings}

The trains running in the networks are homogeneous, each composed of 8 carriages. Each carriage can service 100 passengers or carry 100 units of freight. All trains maintain a constant speed of  300 kilometers per hour. The operational cost per frequency for a line is considered as the unit variable cost of 30,000 CNY per hour multiplied by the duration of the line. The ticket price for passengers and freight is determined by the lengths of the shortest paths in the physical network, and set at 0.7 CNY per kilometer per passenger and 0.2 CNY per kilometer per unit freight. The throughput capacity for each track is capped at 6 trains per hour. The stopping time at each station for each train is fixed to 6 minutes. The transfer times in transfer arcs from a station node to a travel node are uniformly fixed at 30 minutes, while between two consecutive periods, they are set to the duration of the departure period, specifically 4 hours. The duration threshold for passenger demand is set at 1.5 times the shortest travel time, whereas for freight demand it is set at 3 times the shortest travel time. For example, if the shortest travel time from an origin station to a destination station is two hours, then passengers are allowed a maximum travel time of three hours, while freight can take up to six hours to reach the destination.

\section{Computational performance}
We implemented our column generation in Java 8. The RMP in column generation and the MIP  model in price-and-branch are solved using CPLEX version 22.1.1. All experiments are conducted on a single node of a computer cluster. Each node has two 2.6 GHz AMD Rome 7H12 CPU sockets containing 128 CPU cores and 256 GB of RAM. We used 32 cores and 64 GB of RAM for each run. The pricing problems are solved in parallel. 

In the diving heuristic, each column generation phase terminates when there are five consecutive non-improving iterations. The maximum runtime for the MIP solver in the price-and-branch method is set to 30 minutes.

\begin{table}[t]
\renewcommand{\arraystretch}{0.8}
\caption{A comparison of average optimality gap (\%) and solution time (min.) across three ratios for two heuristics on medium and large instances under different demand levels and scenarios}
\centering
\adjustbox{max width=0.9\textwidth}{
\begin{tabular}{cccrrrrrrrrrrr}
\hline
\multirow{3}{*}{\centering \begin{tabular}[c]{@{}c@{}}Demand \\ level\end{tabular}} & \multirow{3}{*}{Scenario} & \multirow{1}{*}{} & \multicolumn{5}{c}{Medium instance}   & \multicolumn{1}{c}{} & \multicolumn{5}{c}{Large instance}                                                                                                           \\ \cline{4-8} \cline{10-14} 
                          &         &              & \multicolumn{2}{c}{Price-and-branch}                      & \multicolumn{1}{c}{} & \multicolumn{2}{c}{Diving heuristic}                      & \multicolumn{1}{c}{} & \multicolumn{2}{c}{Price-and-branch}                      & \multicolumn{1}{c}{} & \multicolumn{2}{c}{Diving heuristic}                      \\ \cline{4-5} \cline{7-8} \cline{10-11} \cline{13-14} 
                          &              &         & \multicolumn{1}{r}{Gap} & \multicolumn{1}{r}{Time} & \multicolumn{1}{c}{} & \multicolumn{1}{r}{Gap} & \multicolumn{1}{r}{Time} & \multicolumn{1}{c}{} & \multicolumn{1}{r}{Gap} & \multicolumn{1}{r}{Time} & \multicolumn{1}{c}{} & \multicolumn{1}{r}{Gap} & \multicolumn{1}{r}{Time} \\ \hline
\multirow{3}{*}{1}        & P             &         & 5.74                       & 30.7                         &                      & 3.84                       & 2.1                          &                      & 23.14                      & 37.3                         &                      & 6.92                       & 18.1                         \\
                          & P+F            &          & 4.83                       & 42.5                         &                      & 3.83                       & 11.8                         &                      & 23.70                      & 106.5                        &                      & 6.97                       & 91.6                         \\
                          & P+F+M            &          & 5.86                       & 171.2                        &                      & 3.75                       & 158.0                        &                      & 64.94                      & 1602.1                       &                      & 6.25                       & 1717.2                       \\ \hline
\multirow{3}{*}{2}        & P               &       & 2.77                       & 30.9                         &                      & 3.19                       & 2.2                          &                      & 14.08                      & 37.9                         &                      & 4.16                       & 21.3                         \\
                          & P+F            &          & 3.22                       & 37.3                         &                      & 3.30                       & 9.8                          &                      & 22.12                      & 89.5                         &                      & 4.61                       & 86.0                         \\
                          & P+F+M           &           & 4.48                       & 145.2                        &                      & 3.04                       & 132.7                        &                      & 48.41                      & 1258.8                       &                      & 4.55                       & 1604.3                       \\ \hline
\multirow{3}{*}{3}        & P              &        & 2.05                       & 31.1                         &                      & 2.50                       & 3.0                          &                      & 9.86                       & 36.4                         &                      & 2.29                       & 23.9                         \\
                          & P+F            &          & 2.25                       & 37.2                         &                      & 2.21                       & 10.1                         &                      & 16.04                      & 82.7                         &                      & 2.62                       & 90.2                         \\
                          & P+F+M          &            & 5.49                       & 141.0                        &                      & 2.17                       & 151.7                        &                      & N/A                     & 1040.9                       &                      & 2.67                       & 1844.0                       \\ \hline
\multirow{3}{*}{4}        & P           &           & 1.71                       & 30.9                         &                      & 2.06                       & 3.5                          &                      & 6.36                       & 35.6                         &                      & 2.16                       & 26.6                         \\
                          & P+F         &             & 1.57                       & 35.9                         &                      & 1.80                       & 10.2                         &                      & 12.66                      & 98.3                         &                      & 2.46                       & 98.4                         \\
                          & P+F+M           &           & 4.16                       & 150.4                        &                      & 1.48                       & 148.3                        &                      & N/A                     & 1100.3                       &                      & 2.60                       & 1931.8                       \\ \hline
\multirow{3}{*}{5}        & P               &       & 1.32                       & 31.0                         &                      & 1.49                       & 3.4                          &                      & 4.55                       & 35.7                         &                      & 1.44                       & 29.3                         \\
                          & P+F            &          & 1.23                       & 35.6                         &                      & 1.37                       & 10.1                         &                      & 15.89                      & 116.2                        &                      & 1.44                       & 118.9                        \\
                          & P+F+M          &            & 2.51                       & 129.1                        &                      & 1.21                       & 107.6                        &                      & N/A                     & 1245.7                       &                      & 1.57                       & 1973.0                       \\ \hline
\multirow{3}{*}{6}        & P            &          & 0.96                       & 31.1                         &                      & 1.17                       & 3.5                          &                      & 11.54                      & 38.7                         &                      & 1.94                       & 34.5                         \\
                          & P+F          &            & 1.15                       & 34.6                         &                      & 1.28                       & 9.0                          &                      & 12.37                      & 106.5                        &                      & 2.20                       & 120.6                        \\
                          & P+F+M        &              & 1.25                       & 126.7                        &                      & 1.28                       & 105.2                        &                      & N/A                     & 1281.5                       &                      & 2.07                       & 2066.2                       \\ \hline
\multirow{3}{*}{7}        & P            &          & 0.52                       & 31.0                         &                      & 0.56                       & 3.6                          &                      & 7.11                       & 38.9                         &                      & 3.13                       & 34.5                         \\
                          & P+F          &            & 0.64                       & 33.6                         &                      & 0.74                       & 7.6                          &                      & 11.40                      & 89.5                         &                      & 1.37                       & 100.8                        \\
                          & P+F+M       &               & 1.17                       & 85.8                         &                      & 0.68                       & 65.2                         &                      & 51.35                      & 1154.5                       &                      & 1.42                       & 1502.5                       \\ \hline
\multirow{3}{*}{8}        & P             &         & 0.47                       & 31.1                         &                      & 0.52                       & 3.5                          &                      & 5.50                       & 42.1                         &                      & 2.31                       & 38.0                         \\
                          & P+F          &            & 0.57                       & 33.7                         &                      & 0.63                       & 7.3                          &                      & 9.30                       & 80.7                         &                      & 2.33                       & 90.1                         \\
                          & P+F+M        &              & 1.21                       & 82.6                         &                      & 0.56                       & 63.5                         &                      & 25.06                      & 778.1                        &                      & 2.41                       & 1169.1                       \\ \hline
\multirow{3}{*}{9}        & P              &        & 0.62                       & 31.0                         &                      & 0.44                       & 3.7                          &                      & 3.80                       & 43.4                         &                      & 2.49                       & 41.6                         \\
                          & P+F           &           & 0.51                       & 32.3                         &                      & 0.48                       & 5.9                          &                      & 7.92                       & 81.9                         &                      & 1.76                       & 91.7                         \\
                          & P+F+M           &           & 1.16                       & 64.2                         &                      & 0.51                       & 42.8                         &                      & 13.56                      & 647.6                        &                      & 1.90                       & 1010.1                       \\ \hline
\multirow{3}{*}{10}       & P              &        & 0.38                       & 31.9                         &                      & 0.33                       & 4.6                          &                      & 4.79                       & 45.9                         &                      & 1.59                       & 53.8                         \\
                          & P+F        &              & 0.39                       & 32.4                         &                      & 0.39                       & 6.0                          &                      & 7.90                       & 80.5                         &                      & 1.62                       & 92.2                         \\
                          & P+F+M        &              & 1.00                       & 53.8                         &                      & 0.42                       & 32.4                         &                      & 12.38                      & 591.7                        &                      & 1.66                       & 900.8                        \\ \hline
\end{tabular}
}
\label{T_AP}
\end{table}

Table \ref{T_AP} presents the optimality gap (\%) and solution time (min.) for the different scenarios and demand levels averaged across three values for the peak-to-off-peak ratio ($R$ of 1.0, 1.5, and 2.0). The gap is computed by comparing the upper bound derived from solving the MP with the solution obtained by the heuristics. In some large instances under scenario P+F+M, the price-and-branch method yields solutions that do not install any lines with a zero profit, indicated with a gap `N/A' (Not Available). 

In all cases, the price-and-branch algorithm fails to achieve optimal solutions within the MIP solver's time limit. As shown in Table \ref{T_AP}, both algorithms can obtain small optimality gaps for medium instances. However, for large instances, the diving algorithm outperforms the price-and-branch algorithm regarding solution quality, albeit requiring more solution time under scenarios P+F and  P+F+M. At each demand level, incorporating more train modes significantly increases the number of integer line variables, as well as the number of nodes and arcs in the change\&go-network, leading to more required diving iterations and, therefore, longer solution time, as well as poorer solution quality in most cases. Additionally, as demand increases, the gap obtained by each algorithm decreases, and the difference between the gaps of the two algorithms for the same instance also diminishes. 

\begin{table}[t]
\renewcommand{\arraystretch}{0.8}
\caption{Average optimality gap (\%) and solution time (min.) across ten demand levels and three ratios  for two heuristics on medium and large instances under different scenarios}
\centering
\adjustbox{max width=.9\textwidth}{
\begin{tabular}{clcrlcrlcrlcr}
\hline
\multirow{3}{*}{Scenario} &  & \multicolumn{5}{c}{Medium instance}                                                                  & \multicolumn{1}{c}{} & \multicolumn{5}{c}{Large instance}                                                                                      \\ \cline{3-7} \cline{9-13} 
                          &  & \multicolumn{2}{c}{Price-and-branch}  & \multicolumn{1}{c}{} & \multicolumn{2}{c}{Diving heuristic}  & \multicolumn{1}{c}{} & \multicolumn{2}{c}{Price-and-branch}                     & \multicolumn{1}{c}{} & \multicolumn{2}{c}{Diving heuristic}  \\ \cline{3-4} \cline{6-7} \cline{9-10} \cline{12-13} 
                          &  & Gap & \multicolumn{1}{r}{Time} & \multicolumn{1}{r}{} & Gap & \multicolumn{1}{r}{Time} & \multicolumn{1}{r}{} & Gap                    & \multicolumn{1}{r}{Time} & \multicolumn{1}{r}{} & Gap & \multicolumn{1}{r}{Time} \\ \hline
P                         &  & 1.65   & 31.07                        &                      & 1.61   & 3.31                         &                      & \multicolumn{1}{r}{9.07}  & 39.20                        &                      & 2.84   & 32.16                        \\
P+F                         &  & 1.64   & 35.52                        &                      & 1.60   & 8.77                         &                      & \multicolumn{1}{r}{13.93} & 93.24                        &                      & 2.74   & 98.03                        \\
P+F+M                        &  & 2.83   & 115.01                       &                      & 1.51   & 100.74                       &                      & \multicolumn{1}{r}{61.57} & 1070.12                      &                      & 2.71   & 1571.92                      \\ \hline
\end{tabular}
}
\label{T_AP2}
\end{table}

Table \ref{T_AP2} shows the average solution time and the gap from Table \ref{T_AP}. The gap of `N/A' is computed by 100\%. As shown in Table \ref{T_AP2}, on average, the diving method outperforms the price-and-branch method regarding solution quality, with this advantage becoming more pronounced for larger instances and scenarios that consider more modes. The diving method demonstrates superior performance across all medium instances regarding solution time. However, for large instances, the average solution time of the diving method significantly increases with the expansion of the line pool size, averaging over 26 hours under scenario P+F+M. Despite the longer time requirement, considering the complexity and strategic nature of the problem, as well as the quality of the solutions obtained, this time is still acceptable. 

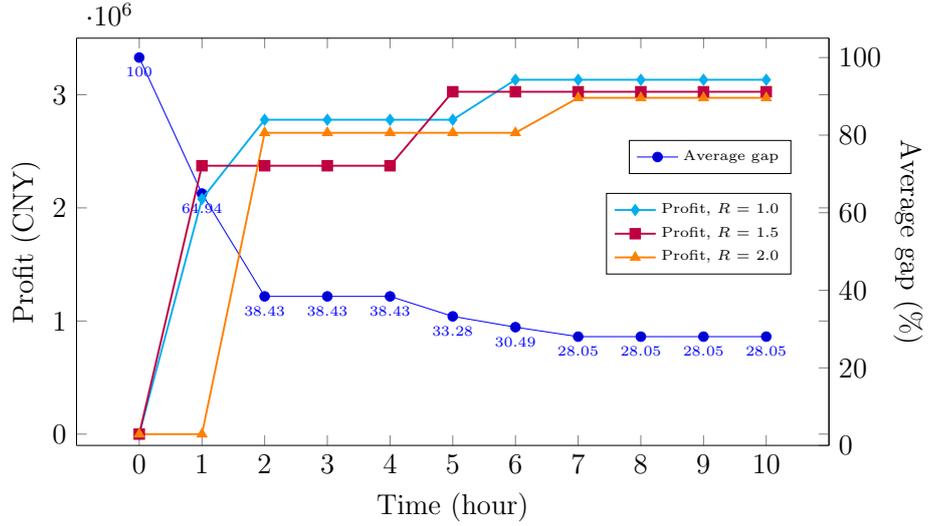
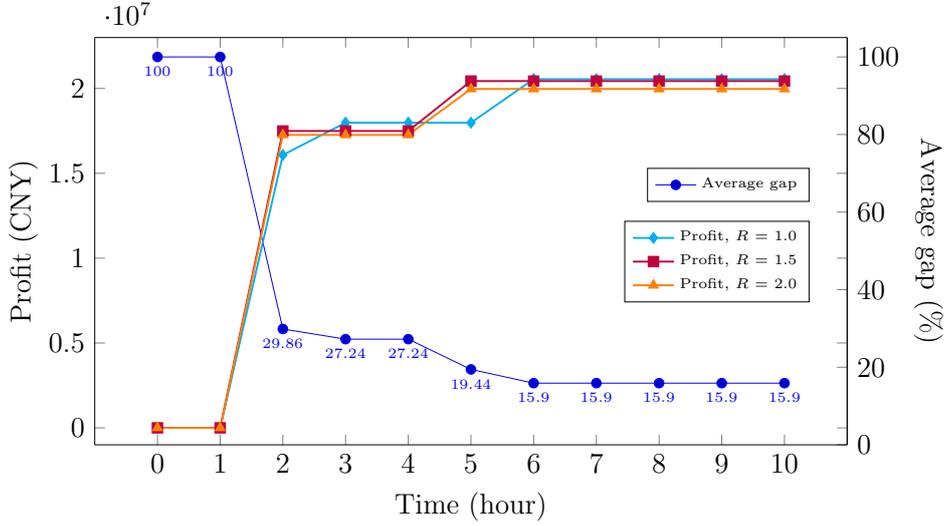
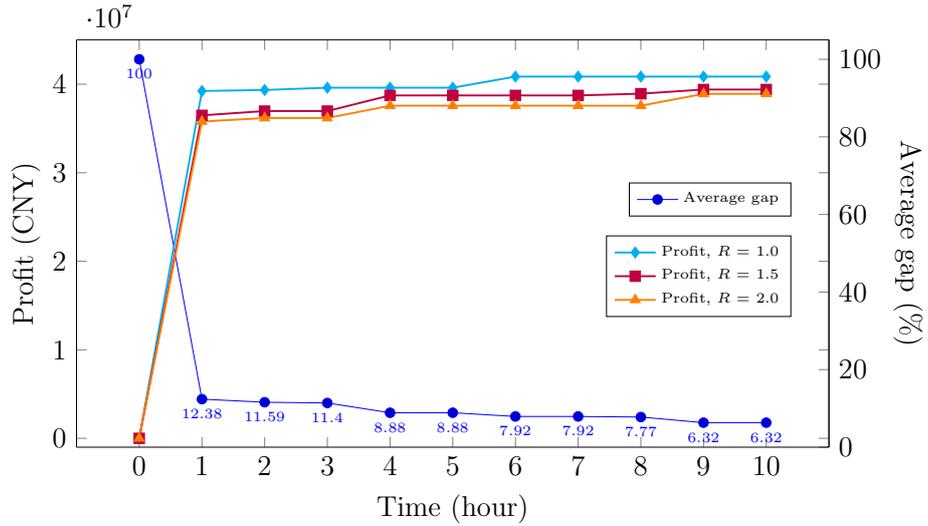
\begin{figure}[htbp]
    \centering
    \begin{subfigure}{0.8\textwidth}
        \centering
         \begin{tikzpicture}[scale=0.9]
        \begin{axis}[
            scale only axis,
             at={(0,0)}, 
            axis y line*=right, 
            height=6cm,
            width=11cm,
            xlabel={Time (hour)},
            ylabel={Profit (CNY)},
            ymin=0, ymax=105,
            xtick=data,
            symbolic x coords={0,1,2,3,4,5,6,7,8,9,10},
            legend style={font=\tiny, at={(0.95,0.75)}, anchor=north east},
            ylabel style={yshift=13cm},
            nodes near coords, 
            every node near coord/.append style={font=\tiny, anchor=north},
        ]
        \addplot
            coordinates {(0,100) (1, 64.94) (2, 38.43) (3,  38.43) (4,  38.43) (5, 33.28) (6, 30.49) (7, 28.05) (8, 28.05) (9, 28.05) (10, 28.05)};
        \addlegendentry{Average gap}
        \end{axis}
        \begin{axis}[
            axis x line=none, 
            scale only axis,
            height=6cm,
            width=11cm,
            ylabel={Average gap (\%)},
            ymin=-100000, ymax=3500000,
            xtick=data,
            symbolic x coords={0,1,2,3,4,5,6,7,8,9,10},
            legend style={font=\tiny, at={(0.95,0.62)}, anchor=north east},
            ylabel style={yshift=-370pt, rotate=180},
        ]
         \addplot[color=cyan, thick, mark=diamond*, mark options={fill=cyan}]
         coordinates  {(0,0) (1, 2078340.7) (2,2778793.85) (3, 2778793.85) (4, 2778793.85) (5,2778793.85) (6, 3132267.49999999) (7, 3132267.49999999) (8, 3132267.49999999) (9, 3132267.49999999) (10, 3132267.49999999)};
        \addlegendentry{Profit, $R$ = 1.0}
        \addplot [color=purple, thick, mark=square*, mark options={fill=purple}]
         coordinates  {(0,0) (1, 2372405.09999999) (2, 2372405.09999999) (3, 2372405.09999999) (4, 2372405.09999999) (5,3026176.19999999) (6, 3026176.19999999) (7, 3026176.19999999) (8, 3026176.19999999) (9, 3026176.19999999) (10, 3026176.19999999)};
        \addlegendentry{Profit, $R$ = 1.5}
        \addplot[color=orange, thick, mark=triangle*, mark options={fill=orange}]
         coordinates  {(0,0) (1, 0) (2, 2663274.4999999) (3, 2663274.4999999) (4, 2663274.4999999) (5, 2663274.4999999) (6, 2663274.4999999) (7, 2972463.85) (8, 2972463.85) (9, 2972463.85) (10, 2972463.85)};
        \addlegendentry{Profit, $R$ = 2.0}
        \end{axis}
    \end{tikzpicture}
    \caption{Large instances, demand level = 1}
    \end{subfigure}
\begin{subfigure}{0.8\textwidth}
        \centering
         \begin{tikzpicture}[scale=0.9]
        \begin{axis}[
            scale only axis,
             at={(0,0)}, 
            axis y line*=right, 
            height=6cm,
            width=11cm,
            xlabel={Time (hour)},
            ylabel={Profit (CNY)},
            ymin=0, ymax=105,
            xtick=data,
            symbolic x coords={0,1,2,3,4,5,6,7,8,9,10},
            legend style={font=\tiny, at={(0.95,0.68)}, anchor=north east},
            ylabel style={yshift=13.3cm},
            nodes near coords, 
            every node near coord/.append style={font=\tiny, anchor=north},
        ]
        \addplot
            coordinates {(0,100) (1, 100) (2, 29.86) (3, 27.24) (4, 27.24) (5, 19.44) (6, 15.9) (7, 15.9) (8, 15.9) (9, 15.9) (10, 15.9)};
        \addlegendentry{Average gap}
        \end{axis}
        \begin{axis}[
            axis x line=none, 
            scale only axis,
            height=6cm,
            width=11cm,
            ylabel={Average gap (\%)},
            ymin=-1000000, ymax=23000000,
            xtick=data,
            symbolic x coords={0,1,2,3,4,5,6,7,8,9,10},
            legend style={font=\tiny, at={(0.95,0.55)}, anchor=north east},
            ylabel style={yshift=-380pt, rotate=180},
        ]
         \addplot[color=cyan, thick, mark=diamond*, mark options={fill=cyan}]
         coordinates  {(0,0) (1, 0) (2, 16079833.6833427) (3, 17981399.3916666) (4, 17981399.3916666) (5, 17981399.3916666) (6, 20548559.7849369) (7, 20548559.7849369) (8, 20548559.7849369) (9, 20548559.7849369) (10, 20548559.7849369)};
        \addlegendentry{Profit, $R$ = 1.0}
        \addplot [color=purple, thick, mark=square*, mark options={fill=purple}]
         coordinates  {(0,0) (1, 0) (2, 17494328.2333333) (3, 17494328.2333333) (4, 17494328.2333333) (5, 20435032.9776166) (6, 20435032.9776166) (7, 20435032.9776166) (8, 20435032.9776166) (9, 20435032.9776166) (10, 20435032.9776166)};
        \addlegendentry{Profit, $R$ = 1.5}
        \addplot[color=orange, thick, mark=triangle*, mark options={fill=orange}]
         coordinates  {(0,0) (1, 0) (2, 17262763.34) (3, 17262763.34) (4, 17262763.34) (5, 19970267.8809824) (6, 19970267.8809824) (7, 19970267.8809824) (8, 19970267.8809824) (9, 19970267.8809824) (10, 19970267.8809824)};
        \addlegendentry{Profit, $R$ = 2.0}
        \end{axis}
    \end{tikzpicture}
    \caption{Large instances, demand level = 5}
    \end{subfigure}

    \begin{subfigure}{0.8\textwidth}
        \centering
         \begin{tikzpicture}[scale=0.9]
        \begin{axis}[
            scale only axis,
             at={(0,0)}, 
            axis y line*=right, 
            height=6cm,
            width=11cm,
            xlabel={Time (hour)},
            ylabel={Profit (CNY)},
            ymin=0, ymax=105,
            xtick=data,
            symbolic x coords={0,1,2,3,4,5,6,7,8,9,10},
            legend style={font=\tiny, at={(0.95,0.65)}, anchor=north east},
            ylabel style={yshift=13cm},
            nodes near coords, 
            every node near coord/.append style={font=\tiny, anchor=north},
        ]
        \addplot
            coordinates {(0,100) (1, 12.38) (2, 11.59) (3, 11.40) (4, 8.88) (5, 8.88) (6, 7.92) (7, 7.92) (8, 7.77) (9, 6.32) (10, 6.32)};
        \addlegendentry{Average gap}
        \end{axis}
        \begin{axis}[
            axis x line=none, 
            scale only axis,
            height=6cm,
            width=11cm,
            ylabel={Average gap (\%)},
            ymin=-1000000, ymax=45000000,
            xtick=data,
            symbolic x coords={0,1,2,3,4,5,6,7,8,9,10},
            legend style={font=\tiny, at={(0.95,0.52)}, anchor=north east},
            ylabel style={yshift=-370pt, rotate=180},
        ]
         \addplot[color=cyan, thick, mark=diamond*, mark options={fill=cyan}]
         coordinates  {(0,0) (1, 39230651.4781464) (2, 39347355.8111666) (3, 39604027.4230052) (4, 39604027.4230052) (5, 39604027.4230052) (6, 40859340.445115) (7, 40859340.445115) (8, 40859340.445115) (9, 40859340.445115) (10, 40859340.445115)};
        \addlegendentry{Profit, $R$ = 1.0}
        \addplot [color=purple, thick, mark=square*, mark options={fill=purple}]
         coordinates  {(0,0) (1, 36480574.237434) (2, 36967122.5194496) (3, 36967122.5194496) (4, 38731293.7952411) (5, 38731293.7952411) (6, 38731293.7952411) (7,38731293.7952411) (8, 38928471.2007268) (9,39399123.378) (10, 39399123.378)};
        \addlegendentry{Profit, $R$ = 1.5}
        \addplot[color=orange, thick, mark=triangle*, mark options={fill=orange}]
         coordinates  {(0,0) (1, 35782211.2099999) (2, 36177469.4961846) (3, 36177469.4961846) (4, 37570836.3193865) (5, 37570836.3193865) (6, 37570836.3193865) (7, 37570836.3193865) (8, 37570836.3193865) (9, 38909715.6912233) (10, 38909715.6912233)};
        \addlegendentry{Profit, $R$ = 2.0}
        \end{axis}
    \end{tikzpicture}
    \caption{Large instances, demand level = 10}
    \end{subfigure}
    \caption{Profit (CNY) obtained by price-and-branch every hour under three ratios and their average optimality gap (\%) for three large instances in scenario P+F+M}
    \label{F_PB}
\end{figure}

Observing substantial differences in the average solution time between price-and-branch and diving methods for large instances under scenario P+F+M, we perform another experiment where we extend the solution time limit in the branching step of price-and-branch to 10 hours and conduct further tests on large instances with demand level 1, 5, and 10. Feasible solutions obtained by the solver are recorded every hour, and the average gaps are calculated, as depicted in Figure \ref{F_PB}. Feasible solutions are obtained for all cases when the MIP solver's runtime reaches 2 hours. As the solution time increases, the average gap decreases at a diminishing rate and eventually stabilizes until the time limit is reached. The average gaps for large instances with demand levels 1, 5, and 10 finally converge to 28.05\%, 15.90\%, and 6.43\%, respectively, compared to 6.25\%, 1.57\%, and 1.66\% achieved by the diving method. While the solution quality of price-and-branch improves with extended solution time, it remains inferior to that of the diving method.

\begin{figure}[htbp]
    \centering
    \begin{subfigure}{0.45\textwidth}
        \centering
         \begin{tikzpicture}[scale=1]
        \begin{axis}[
            ybar, 
            bar width=10pt,
            scale only axis,
            height=5cm,
            width=6cm,
            at={(0,0)}, 
            axis y line*=right, 
            xlabel={Demand level},
            ylabel={Solution time (min.)},
            ymin=1, ymax=5.3, 
            xtick=data,
            symbolic x coords={1,2,3,4,5,6,7,8,9,10},
            legend style={font=\tiny, at={(0.55,0.98)}, anchor=north west},
            ylabel style={yshift=5pt, rotate=180},
        ]
        \addplot
            coordinates {(1,2.09) (2,2.24) (3,3.05) (4, 3.48) (5,3.42) (6,3.54) (7, 3.57) (8,3.47) (9,3.73) (10,4.55)};
        \addlegendentry{Solution time}
        \end{axis}
        \begin{axis}[
            axis x line=none, 
            scale only axis,
            height=5cm,
            width=6cm,
            xlabel={Demand level},
            ylabel={Number of paths},
            ymin=0, ymax=22000,
            xtick=data,
            symbolic x coords={1,2,3,4,5,6,7,8,9,10},
            legend style={font=\tiny, at={(0.05,0.97)}, anchor=north west},
            ylabel style={yshift=-5pt},
        ]
        \addplot[color=red, mark=square*]
            coordinates {(1,5843) (2,7289) (3,9899) (4, 10442) (5,11221) (6,11909) (7, 12604) (8,13897) (9,14763) (10,17887)};
        \addlegendentry{Passenger path}
        \end{axis}
    \end{tikzpicture}
    \caption{Medium instances, scenario P}
    \end{subfigure}
    \hspace{0.05\textwidth}
    \begin{subfigure}{0.45\textwidth}
        \centering
        \begin{tikzpicture}[scale=1]
        \begin{axis}[
            axis x line=none, 
            ybar, 
            bar width=10pt,
            scale only axis,
            height=5cm,
            width=6cm,
            at={(0,0)}, 
            axis y line*=right, 
            ylabel={Solution time (min.)},
            ymin=0, ymax=65, 
            xtick=data,
            symbolic x coords={1,2,3,4,5,6,7,8,9,10},
            legend style={font=\tiny, at={(0.55,0.98)}, anchor=north west},
            ylabel style={yshift=5pt, rotate=180},
        ]
        \addplot
            coordinates {(1,18.06) (2,21.32) (3,23.91) (4, 26.61) (5,29.30) (6,34.5) (7, 34.46) (8,37.99) (9,41.64) (10,53.76)};
        \addlegendentry{Solution time}
        
        \end{axis}
        \begin{axis}[
            scale only axis,
            height=5cm,
            width=6cm,
            xlabel={Demand level},
            ylabel={Number of paths},
            ymin=10000, ymax=55000,
            xtick=data,
            symbolic x coords={1,2,3,4,5,6,7,8,9,10},
            legend style={font=\tiny, at={(0.05,0.97)}, anchor=north west},
            ylabel style={yshift=-5pt},
        ]
        \addplot[color=red, mark=square*]
            coordinates {(1,14716) (2,21597) (3,24831) (4, 27855) (5,31662) (6,34282) (7, 37268) (8,42677) (9,43289) (10,46652)};
        \addlegendentry{Passenger path}
        \end{axis}
    \end{tikzpicture}
    \caption{Large instances, scenario P}
    \end{subfigure}
    \begin{subfigure}{0.45\textwidth}
        \centering
         \begin{tikzpicture}[scale=1]
        \begin{axis}[
            axis x line=none, 
            ybar, 
            bar width=10pt,
            scale only axis,
            width=0.5\textwidth,
            height=5cm,
            width=6cm,
            at={(0,0)}, 
            axis y line*=right, 
            ylabel={Solution time (min.)},
            ymin=0, ymax=15, 
            xtick=data,
            symbolic x coords={1,2,3,4,5,6,7,8,9,10},
            legend style={font=\tiny, at={(0.55,0.96)}, anchor=north west},
            ylabel style={yshift=5pt, rotate=180},
        ]
        \addplot
            coordinates {(1,11.77) (2,9.79) (3,10.08) (4, 10.16) (5,10.07) (6,9.00) (7, 7.63) (8,7.27) (9,5.91) (10,6.01)};
        \addlegendentry{Solution time}
        
        \end{axis}
        \begin{axis}[
            scale only axis,
            width=0.5\textwidth,
            height=5cm,
            width=6cm,
            xlabel={Demand level},
            ylabel={Number of paths},
            ymin=5000, ymax=22000,
            xtick=data,
            symbolic x coords={1,2,3,4,5,6,7,8,9,10},
            legend style={font=\tiny, at={(0.05,0.98)}, anchor=north west},
            ylabel style={yshift=-5pt},
        ]
        \addplot[color=red, mark=square*]
            coordinates {(1,6684) (2,7854) (3,10631) (4, 10675) (5,11576) (6,12308) (7, 12660) (8,14513) (9,15216) (10,18117)};
        \addlegendentry{Passenger path}
        \addplot[color=blue, mark=*]
            coordinates {(1,11840) (2,10259) (3,11301) (4,11944) (5,11955) (6,11600) (7,10347) (8,10092) (9,8002) (10,7755)};
        \addlegendentry{Freight path}
        \end{axis}
    \end{tikzpicture}
    \caption{Medium instances, scenario P+F}
    \end{subfigure}
    \hspace{0.05\textwidth}
    \begin{subfigure}{0.45\textwidth}
        \centering
         \begin{tikzpicture}[scale=1]
        \begin{axis}[
            axis x line=none, 
            ybar, 
            bar width=10pt,
            scale only axis,
            width=0.5\textwidth,
            height=5cm,
            width=6cm,
            at={(0,0)}, 
            axis y line*=right, 
            ylabel={Solution time (min.)},
            ymin=25, ymax=150, 
            xtick=data,
            symbolic x coords={1,2,3,4,5,6,7,8,9,10},
            legend style={font=\tiny, at={(0.55,0.96)}, anchor=north west},
            ylabel style={yshift=5pt, rotate=180},
        ]
        \addplot
            coordinates {(1,91.57) (2,85.98) (3,90.23) (4, 98.39) (5,118.89) (6,120.56) (7, 100.76) (8,90.06) (9,91.69) (10,92.21)};
        \addlegendentry{Solution time}
        
        \end{axis}
        \begin{axis}[
            scale only axis,
            width=0.5\textwidth,
            height=5cm,
            width=6cm,
            xlabel={Demand level},
            ylabel={Number of paths},
            ymin=10000, ymax=60000,
            xtick=data,
            symbolic x coords={1,2,3,4,5,6,7,8,9,10},
            legend style={font=\tiny, at={(0.05,0.98)}, anchor=north west},
            ylabel style={yshift=-5pt},
        ]
        \addplot[color=red, mark=square*]
            coordinates {(1,15434) (2,21473) (3,25193) (4, 29089) (5,31531) (6,35039) (7, 38346) (8,40335) (9,43909) (10,45384)};
        \addlegendentry{Passenger path}
        \addplot[color=blue, mark=*]
            coordinates {(1,32343) (2,31601) (3,35229) (4,38021) (5,44775) (6,42223) (7,36725) (8,34882) (9,33101) (10,32641)};
        \addlegendentry{Freight path}
        \end{axis}
    \end{tikzpicture}
    \caption{Large instances, scenario P+F}
    \end{subfigure}
     \begin{subfigure}{0.45\textwidth}
        \centering
         \begin{tikzpicture}[scale=1]
        \begin{axis}[
            axis x line=none, 
            ybar, 
            bar width=10pt,
            scale only axis,
            width=0.5\textwidth,
            height=5cm,
            width=6cm,
            at={(0,0)}, 
            axis y line*=right, 
            ylabel={Solution time (min.)},
            ymin=0, ymax=200, 
            xtick=data,
            symbolic x coords={1,2,3,4,5,6,7,8,9,10},
            legend style={font=\tiny, at={(0.55,0.96)}, anchor=north west},
            ylabel style={yshift=5pt, rotate=180},
        ]
        \addplot
            coordinates {(1,158) (2,132.72) (3,151.7) (4, 148.31) (5,107.64) (6,105.18) (7,65.21) (8,63.52) (9,42.77) (10,32.38)};
        \addlegendentry{Solution time}
        
        \end{axis}
        \begin{axis}[
            scale only axis,
            width=0.5\textwidth,
            height=5cm,
            width=6cm,
            xlabel={Demand level},
            ylabel={Number of paths},
            ymin=15000, ymax=45000,
            xtick=data,
            symbolic x coords={1,2,3,4,5,6,7,8,9,10},
            legend style={font=\tiny, at={(0.05,0.98)}, anchor=north west},
            ylabel style={yshift=-5pt},
        ]
        \addplot[color=red, mark=square*]
            coordinates {(1,17278) (2,17459) (3,22669) (4, 24210) (5,24191) (6,26517) (7, 26823) (8,29762) (9,30634) (10,38148)};
        \addlegendentry{Passenger path}
        \addplot[color=blue, mark=*]
            coordinates {(1,35235) (2,28889) (3,29827) (4,33141) (5,29247) (6,30136) (7,23917) (8,24423) (9,19481) (10,17010)};
        \addlegendentry{Freight path}
        \end{axis}
    \end{tikzpicture}
    \caption{Medium instances, scenario P+F+M}
    \end{subfigure}
    \hspace{0.05\textwidth}
    \begin{subfigure}{0.45\textwidth}
        \centering
           \begin{tikzpicture}[scale=1]
        \begin{axis}[
            axis x line=none, 
            ybar, 
            bar width=10pt,
            scale only axis,
            width=0.5\textwidth,
            height=5cm,
            width=6cm,
            at={(0,0)}, 
            axis y line*=right, 
            ylabel={Solution time (min.)},
            ymin=500, ymax=2500, 
            xtick=data,
            symbolic x coords={1,2,3,4,5,6,7,8,9,10},
            legend style={font=\tiny, at={(0.55,0.96)}, anchor=north west},
            ylabel style={yshift=5pt, rotate=180},
        ]
        \addplot
            coordinates {(1,1717.24) (2,1604.33) (3,1844.01) (4, 1931.76) (5,1973.04) (6,2066.23) (7, 1502.52) (8,1169.14) (9,1010.11) (10,900.78)};
        \addlegendentry{Solution time}
        
        \end{axis}
        \begin{axis}[
            scale only axis,
            width=0.5\textwidth,
            height=5cm,
            width=6cm,
            xlabel={Demand level},
            ylabel={Number of paths},
            ymin=30000, ymax=130000,
            xtick=data,
            symbolic x coords={1,2,3,4,5,6,7,8,9,10},
            legend style={font=\tiny, at={(0.05,0.98)}, anchor=north west},
            ylabel style={yshift=-5pt},
        ]
        \addplot[color=red, mark=square*]
            coordinates {(1,37226) (2,47750) (3,49726) (4, 54565) (5,60342) (6,63001) (7, 72460) (8,77896) (9,83080) (10,88371)};
        \addlegendentry{Passenger path}
        \addplot[color=blue, mark=*]
            coordinates {(1,84101) (2,79934) (3,86878) (4,93381) (5,103345) (6,99395) (7,93510) (8,89196) (9,80558) (10,79814)};
        \addlegendentry{Freight path}
        \end{axis}
    \end{tikzpicture}
    \caption{Large instances, scenario P+F+M}
    \end{subfigure}
    \caption{Average number of generated paths and solution time (min.) across three ratios on medium and large instances for different demand levels and scenarios}
    \label{Diving}
\end{figure}
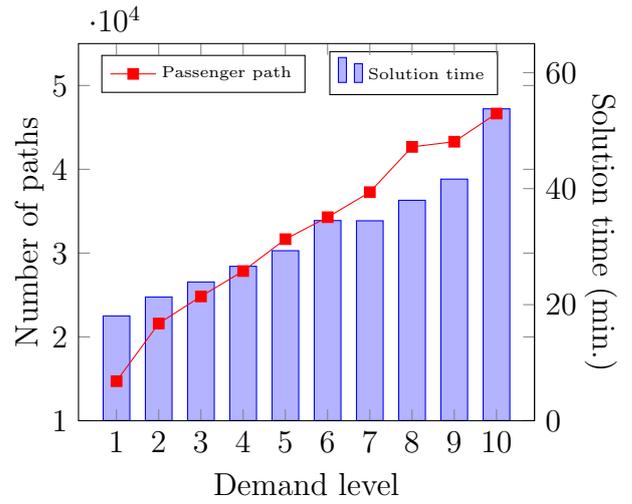
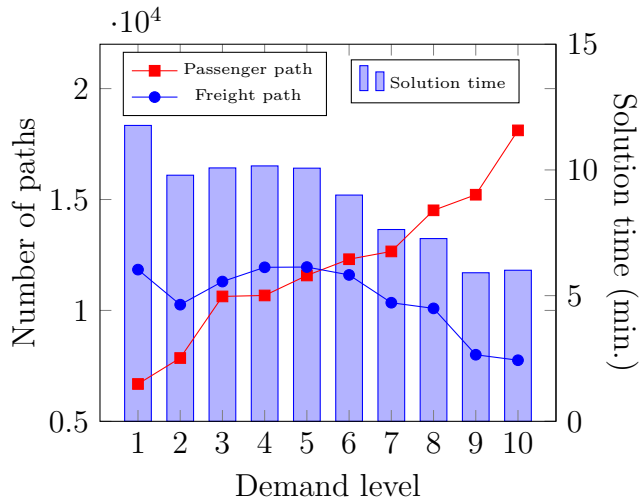
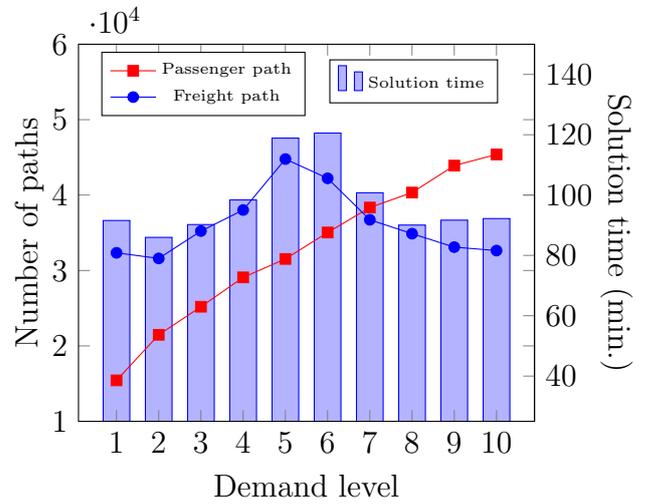
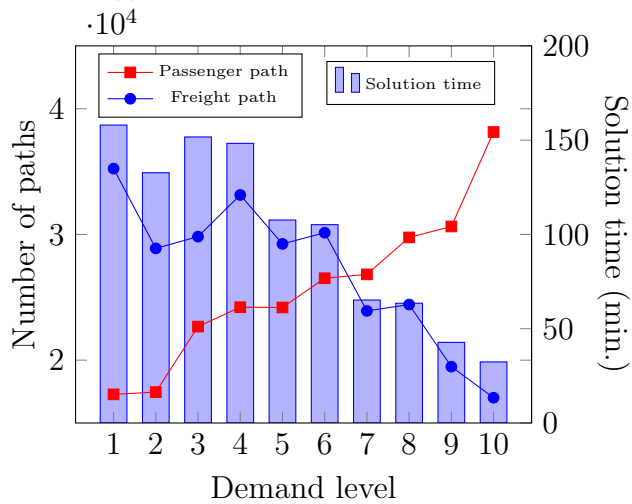
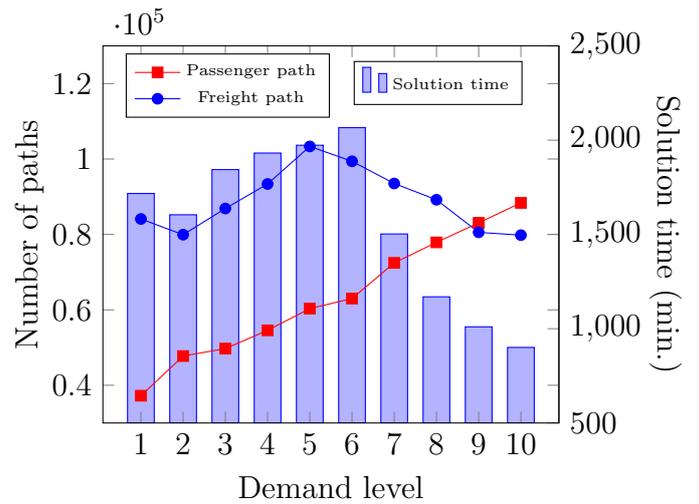

To further analyze what drives the solution time of the diving heuristic, Figure \ref{Diving} illustrates the number of generated passenger and freight paths during the solving process. There appears to be a strong correlation between the solution time and the number of generated paths, especially for the freight paths. In passenger-only scenarios, the solution time and the number of generated paths gradually increase with rising demand levels. However, in combined passenger-freight scenarios, the solution time initially increases and then decreases as demand levels rise, mirroring the trend in the number of generated freight paths. Freight paths can start or end in any period, requiring a search across the entire network considering interconnected period layers, whereas passenger paths are searched within specific period layers. This makes the freight pricing problem more time-consuming and explains why the trend in solution time aligns with the trend in the number of generated freight paths. 

Additionally, in high passenger demand, railway lines are primarily allocated to passenger service, resulting in limited capacity for freight transportation and reduced potential for generating freight paths. Conversely, surplus line capacity allows for more freight path generation in instances with lower passenger demand. This explains why the number of freight paths initially increases and then decreases as demand levels rise.

Based on the results presented and the trade-off between solution time and solution quality, we determine that the diving heuristic performs better. Therefore, we will proceed with this algorithm for further analysis.

\section{Management insights}

\begin{table}[t]
\centering
\renewcommand{\arraystretch}{0.8}
\caption{Average service level (SL, \%) and travel time (TT, hr) across ten demand levels under different ratios and scenarios}
\adjustbox{max width=1\textwidth}{
\begin{tabular}{ccccccccccccccccccccccc}
\hline
\multirow{3}{*}{Demand}      &  & \multirow{3}{*}{Scenario} &  & \multicolumn{19}{c}{peak-to-off-peak ratio of passenger demand $R$}                                                                                              \\ \cline{5-23} 
                           &  &                           &  & \multicolumn{3}{c}{1.0} &  & \multicolumn{3}{c}{1.5} &  & \multicolumn{3}{c}{2.0} &  & \multicolumn{3}{c}{2.5} &  & \multicolumn{3}{c}{3.0} \\ \cline{5-7} \cline{9-11} \cline{13-15} \cline{17-19} \cline{21-23} 
                           &  &                           &  & SL       &     & TT     &  & SL       &     & TT     &  & SL       &     & TT     &  & SL       &     & TT     &  & SL       &     & TT     \\ \hline
\multirow{3}{*}{Passenger} &  & P                         &  & 88.85    &     & 1.30   &  & 88.51    &     & 1.30   &  & 87.38    &     & 1.30   &  & 86.86    &     & 1.30   &  & 85.76    &     & 1.31   \\
                           &  & P+F                       &  & 88.91    &     & 1.30   &  & 88.55    &     & 1.30   &  & 87.37    &     & 1.31   &  & 86.89    &     & 1.31   &  & 85.83    &     & 1.31   \\
                           &  & P+F+M                     &  & 89.11    &     & 1.31   &  & 88.79    &     & 1.30   &  & 88.46    &     & 1.31   &  & 87.38    &     & 1.31   &  & 86.14    &     & 1.31   \\ \hline
\multirow{2}{*}{Freight}   &  & P+F                       &  & 24.74    &     & 1.22   &  & 23.46    &     & 1.19   &  & 22.57    &     & 1.19   &  & 20.99    &     & 1.17   &  & 20.81    &     & 1.14   \\
                           &  & P+F+M                     &  & 25.03    &     & 1.35   &  & 24.64    &     & 1.29   &  & 23.08    &     & 1.27   &  & 22.34    &     & 1.24   &  & 22.25    &     & 1.21   \\ \hline
\end{tabular}
}
\label{T_TT}
\end{table}

In this section, we use the developed diving heuristic to address the following questions: 

\begin{enumerate}[i)]
    \item What benefits can the railway passenger transportation system derive from integrating freight services? (Section 7.1)
    \item Is multi-period line planning necessary to accommodate time-varying passenger demand in an integrated passenger-freight railway system? (Section 7.2)
    \item What is the impact of different passenger demand patterns on the characteristics of the line plan? (Section 7.3)
\end{enumerate}

\subsection{Benefits of integrating freight}
We firstly present the average service level (SL, measured in \%), defined as the proportion of demand that is served, and average travel time (TT, measured in hours) for both passengers and freight under different 
peak-to-off-peak ratios and scenarios in Table \ref{T_TT}. Each entry in the table corresponds to the average outcome of ten demand levels. 

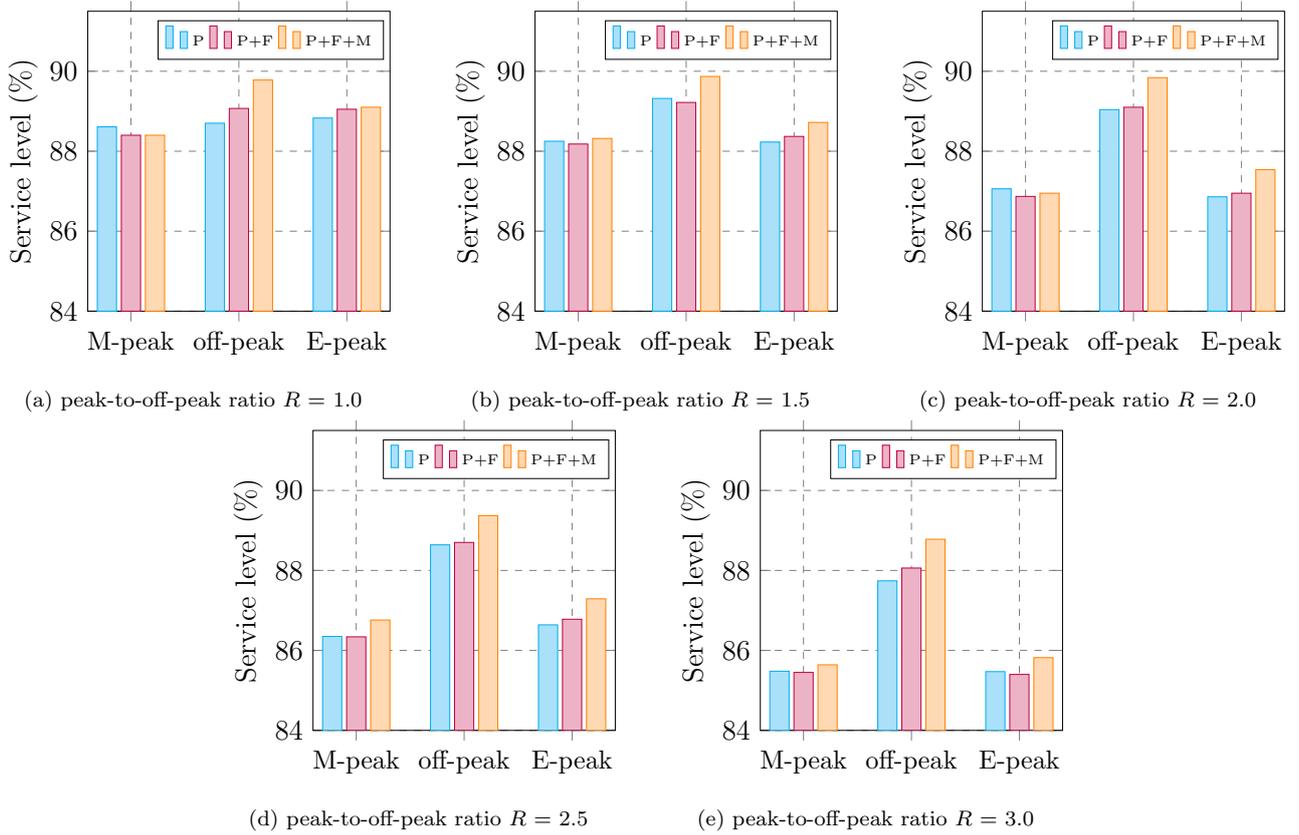
\begin{figure}[t]
    \centering
    \begin{subfigure}{0.3\textwidth} 
        \centering
       \begin{tikzpicture}[scale=0.9]
            \begin{axis}[
                ybar,
                symbolic x coords={M-peak, off-peak, E-peak},
                ylabel={Service level (\%)},
                ymin=84, ymax=91.5,
                xtick=data,
                xticklabel style={font=\small},
                nodes near coords={}, 
                nodes near coords align={vertical},
                enlarge x limits={0.2, 0.2},
                ylabel style={yshift=-3pt},
                legend style={font=\tiny, at={(0.98,0.97)}, anchor=north east, legend columns=3},
                grid=major,
                grid style={dashed, gray},
                bar width=8pt,
                width=6cm,
                height=6cm
            ]
                \addplot[color=cyan, fill=cyan!30] coordinates {(M-peak,88.61) (off-peak,88.70) (E-peak,88.83)};
                \addplot[color=purple, fill=purple!30] coordinates {(M-peak,88.40) (off-peak,89.07) (E-peak,89.05)};
                \addplot[color=orange, fill=orange!30] coordinates {(M-peak,88.40) (off-peak,89.78) (E-peak,89.10)};
                \legend{P, P+F, P+F+M}
            \end{axis}
        \end{tikzpicture}
        \caption{ peak-to-off-peak ratio $R$ = 1.0}
    \end{subfigure}
    \hspace{0.03\textwidth}
    \begin{subfigure}{0.3\textwidth} 
        \centering
       \begin{tikzpicture}[scale=0.9]
            \begin{axis}[
                ybar,
                symbolic x coords={M-peak, off-peak, E-peak},
                ymin=84, ymax=91.5,
                xtick=data,
                xticklabel style={font=\small},
                nodes near coords={}, 
                nodes near coords align={vertical},
                enlarge x limits={0.2, 0.2},
                ylabel={Service level (\%)},
                ylabel style={yshift=-3pt},
                legend style={font=\tiny, at={(0.98,0.97)}, anchor=north east, legend columns=3},
                grid=major,
                grid style={dashed, gray},
                 bar width=8pt,
                width=6cm,
                height=6cm
            ]
                \addplot[color=cyan, fill=cyan!30] coordinates {(M-peak,88.25) (off-peak,89.32) (E-peak,88.23)};
                \addplot[color=purple, fill=purple!30] coordinates {(M-peak,88.18) (off-peak,89.22) (E-peak,88.37)};
                \addplot[color=orange, fill=orange!30] coordinates {(M-peak,88.32) (off-peak,89.87) (E-peak,88.72)};
                \legend{P, P+F, P+F+M}
            \end{axis}
        \end{tikzpicture}
        \caption{peak-to-off-peak ratio $R$ = 1.5}
    \end{subfigure}
    \hspace{0.03\textwidth}
    \begin{subfigure}{0.3\textwidth} 
        \centering
       \begin{tikzpicture}[scale=0.9]
            \begin{axis}[
                ybar,
                symbolic x coords={M-peak, off-peak, E-peak},
                ymin=84, ymax=91.5,
                xtick=data,
                xticklabel style={font=\small},
                nodes near coords={}, 
                nodes near coords align={vertical},
                enlarge x limits={0.2, 0.2},
                ylabel={Service level (\%)},
                ylabel style={yshift=-3pt},
                legend style={font=\tiny, at={(0.98,0.97)}, anchor=north east, legend columns=3},
                grid=major,
                grid style={dashed, gray},
                 bar width=8pt,
                width=6cm,
                height=6cm
            ]
                \addplot[color=cyan, fill=cyan!30] coordinates {(M-peak,87.06) (off-peak,89.04) (E-peak,86.86)};
                \addplot[color=purple, fill=purple!30] coordinates {(M-peak,86.87) (off-peak,89.10) (E-peak,86.95)};
                \addplot[color=orange, fill=orange!30] coordinates {(M-peak,86.95) (off-peak,89.84) (E-peak,87.54)};
                \legend{P, P+F, P+F+M}
            \end{axis}
        \end{tikzpicture}
        \caption{peak-to-off-peak ratio $R$ = 2.0}
    \end{subfigure}
    \newline
    \begin{subfigure}{0.3\textwidth} 
        \centering
       \begin{tikzpicture}[scale=0.9]
            \begin{axis}[
                ybar,
                symbolic x coords={M-peak, off-peak, E-peak},
                ymin=84, ymax=91.5,
                xtick=data,
                xticklabel style={font=\small},
                nodes near coords={}, 
                nodes near coords align={vertical},
                enlarge x limits={0.2, 0.2},
                ylabel={Service level (\%)},
                ylabel style={yshift=-3pt},
                legend style={font=\tiny, at={(0.98,0.97)}, anchor=north east, legend columns=3},
                grid=major,
                grid style={dashed, gray},
                bar width=8pt,
                width=6cm,
                height=6cm
            ]
                \addplot[color=cyan, fill=cyan!30] coordinates {(M-peak,86.35) (off-peak,88.64) (E-peak,86.64)};
                \addplot[color=purple, fill=purple!30] coordinates {(M-peak,86.34) (off-peak,88.70) (E-peak,86.78)};
                \addplot[color=orange, fill=orange!30] coordinates {(M-peak,86.76) (off-peak,89.37) (E-peak,87.29)};
                \legend{P, P+F, P+F+M}
            \end{axis}
        \end{tikzpicture}
        \caption{peak-to-off-peak ratio $R$ = 2.5}
    \end{subfigure}
     \hspace{0.03\textwidth}
    \begin{subfigure}{0.3\textwidth} 
        \centering
       \begin{tikzpicture}[scale=0.9]
            \begin{axis}[
                ybar,
                symbolic x coords={M-peak, off-peak, E-peak},
                ymin=84, ymax=91.5,
                xtick=data,
                xticklabel style={font=\small},
                nodes near coords={}, 
                nodes near coords align={vertical},
                enlarge x limits={0.2, 0.2},
                legend style={font=\tiny, at={(0.98,0.97)}, anchor=north east, legend columns=3},
                ylabel={Service level (\%)},
                ylabel style={yshift=-3pt},
                grid=major,
                grid style={dashed, gray},
                bar width=8pt,
                width=6cm,
                height=6cm
            ]
               \addplot[color=cyan, fill=cyan!30] coordinates {(M-peak,85.48) (off-peak,87.74) (E-peak,85.47)};
                \addplot[color=purple, fill=purple!30] coordinates {(M-peak,85.45) (off-peak,88.06) (E-peak,85.40)};
                \addplot[color=orange, fill=orange!30] coordinates {(M-peak,85.64) (off-peak,88.78) (E-peak,85.82)};
                \legend{P, P+F, P+F+M}
            \end{axis}
        \end{tikzpicture}
        \caption{peak-to-off-peak ratio $R$ = 3.0}
    \end{subfigure}

    \caption{Average passenger service level (\%) in each period across ten demand levels under different ratios and scenarios}
    \label{F_Ser2}
\end{figure} 

As shown in Table \ref{T_TT}, integrating freight does not compromise passenger service levels, as there is no decrease in passenger service levels when the railway system is used to transport passengers and freight. In scenarios with trains that carry both passengers and freight, there can even be a slight improvement in passenger service because the additional revenue from freight can be used to operate extra trains. At the same time, the average travel time for passengers remains nearly unchanged, which can be explained by the fact that our model does not explicitly minimize travel time but enforces a maximum allowed detour using constraints. 

 Table \ref{T_TT} also shows that passenger service levels decrease when the peak-to-off-peak ratio increases, which suggests that time-varying passenger demand negatively impacts service levels. To analyze this further, Figure \ref{F_Ser2} illustrates the average passenger service levels in the different periods. Overall, passenger service levels decline within both peak and off-peak periods as the demand becomes less and less uniform. This decline is primarily due to the conflict between increased demand and saturated capacity during the M-peak and E-peak periods. In off-peak periods, the reduced demand poses challenges in pooling sufficient demand to operate trains profitably. However, the introduction of mixed trains allows for flexible combinations of passenger and freight demand, significantly improving passenger service levels during off-peak periods and contributing to overall higher passenger service levels under scenario P+F+M.

In terms of freight service, as depicted in Table \ref{T_TT}, scenario P+F+M shows higher levels of freight service compared to scenario P+F. However, the average travel time for freight is longer in scenario P+F+M than in scenario P+F, indicating that offering mixed trains can improve freight service levels at the expense of a longer average travel duration. Additionally, as the peak-to-off-peak ratio increases, the system tends to prioritize shorter-duration freight services, reducing the average transportation time for freight.

\begin{figure}[t]
    \centering
    \begin{subfigure}{0.45\textwidth}
        \centering
         \begin{tikzpicture}[scale=1]
        \begin{axis}[
            xlabel={peak-to-off-peak ratio $R$},
            ylabel={Utilization (\%)},
            symbolic x coords={1.0,1.5,2.0,2.5,3.0},
            ymin=89, ymax=93,
            xtick=data,
            enlarge x limits={0.15},
            legend style={font=\tiny, at={(0.05,0.98)}, anchor=north west, legend columns=1},
            ylabel style={yshift=-1mm},
            xlabel style={yshift=0mm},
            grid=major,
            grid style={dashed, gray},
            every axis plot/.append style={thick},
            tick label style={font=\small},
            width=7.5cm,
            height=7cm
        ]

            \addplot[color=cyan, mark=diamond*, mark options={fill=cyan}] coordinates {
                 (1.0,89.8) (1.5,89.92) (2.0,90.6) (2.5,90.69) (3.0,91.1)
            };
            \addplot[color=purple, mark=square*, mark options={fill=purple}] coordinates {
                 (1.0,90.05) (1.5,90.49) (2.0,90.85) (2.5,90.99) (3.0,91.46)
            };
            \addplot[color=orange, mark=triangle*, mark options={fill=orange}] coordinates {
                 (1.0,90.73) (1.5,91.42) (2.0,91.74) (2.5,92.00) (3.0,92.19)
            };
            \legend{P, P+F, P+F+M}
        \end{axis}
    \end{tikzpicture}
    \caption{Passenger total utilization}
    \label{F_Util 1}
    \end{subfigure}
    \hspace{0.05\textwidth}
    \begin{subfigure}{0.45\textwidth}
        \centering
        \begin{tikzpicture}[scale=1]
        \begin{axis}[
            xlabel={peak-to-off-peak ratio $R$},
            ylabel={Utilization (\%)},
            symbolic x coords={1.0,1.5,2.0,2.5,3.0},
            ymin=96, ymax=100,
            xtick=data,
            enlarge x limits={0.15},
            legend style={font=\tiny, at={(0.05,0.98)}, anchor=north west, legend columns=1},
            ylabel style={yshift=-1mm},
            xlabel style={yshift=0mm},
            grid=major,
            grid style={dashed, gray},
            every axis plot/.append style={thick},
            tick label style={font=\small},
            width=7.5cm,
            height=7cm
        ]
        
            \addplot[color=purple, mark=square*, mark options={fill=purple}] coordinates {
                (1.0,98.16) (1.5,97.05) (2.0,97.16) (2.5,96.81) (3.0,97.46)
            };
            \addplot[color=orange, mark=triangle*, mark options={fill=orange}] coordinates {
                (1.0,98.67) (1.5,98.30) (2.0,98.18) (2.5,98.40) (3.0,97.94)
            };
            \legend{P+F, P+F+M}
        \end{axis}
    \end{tikzpicture}
    \caption{Freight total utilization}
    \label{F_Util 2}
    \end{subfigure}
    \caption{Average total utilization (\%) across ten demand levels under different ratios and scenarios}
    \label{F_Util}
\end{figure}
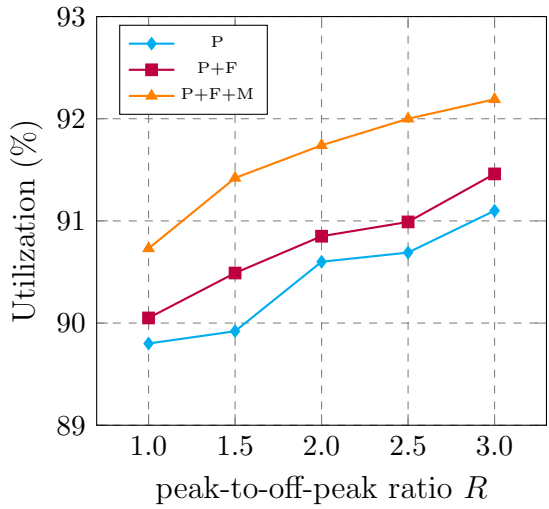
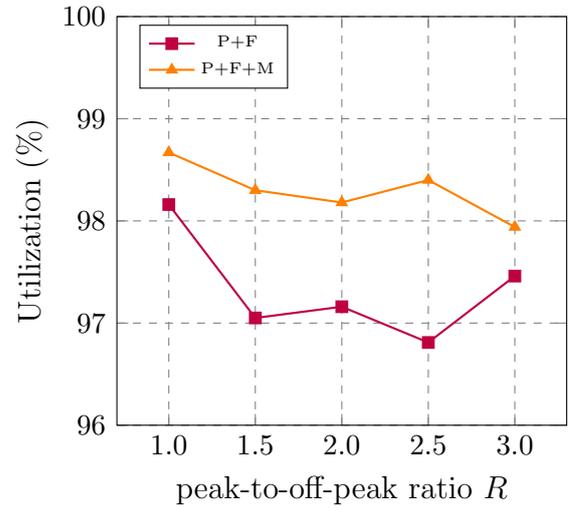

Next, we present the average passenger and freight utilization across ten demand levels under different scenarios and peak-to-off-peak ratios in Figure \ref{F_Util}, with \ref{F_Util 1} depicting passenger utilization and \ref{F_Util 2} depicting freight utilization. Passenger utilization is defined as the average loading rate of travel arcs related to installed lines with positive passenger capacity. Similarly, freight utilization is defined as the average loading rate of travel arcs associated with installed lines with positive freight capacity.

\begin{table}[t]
\renewcommand{\arraystretch}{0.8}
\centering
\caption{Average passenger utilization (\%) in each period across ten demand levels under different ratios and scenarios}
\adjustbox{max width=1\textwidth}{
\begin{tabular}{ccccccccccccccccccc}
\hline
\multirow{2}{*}{$R$} & \multirow{2}{*}{} & \multicolumn{5}{c}{Scenario P}       &  & \multicolumn{5}{c}{Scenario P+F}       &  & \multicolumn{5}{c}{Scenario P+F+M}       \\ \cline{3-7} \cline{9-13} \cline{15-19} 
                       &                   & M-peak &  & Off-peak &  & E-peak &  & M-peak &  & Off-peak &  & E-peak &  & M-peak &  & Off-peak &  & E-peak \\ \hline
1.0                      &                   & 89.82  &  & 89.79    &  & 89.89  &  & 91.16  &  & 89.99    &  & 90.63  &  & 91.30  &  & 90.03    &  & 90.66  \\
1.5                    &                   & 90.11  &  & 89.90    &  & 89.87  &  & 91.63  &  & 90.59    &  & 91.26  &  & 91.98  &  & 90.88    &  & 90.76  \\
2.0                      &                   & 90.78  &  & 89.85    &  & 91.23  &  & 91.51  &  & 91.02    &  & 91.49  &  & 92.45  &  & 91.26    &  & 91.61  \\
2.5                    &                   & 91.52  &  & 89.44    &  & 91.32  &  & 92.13  &  & 89.15    &  & 91.80  &  & 93.16  &  & 89.81    &  & 91.90  \\
3.0                      &                   & 91.92  &  & 89.69    &  & 91.58  &  & 92.73  &  & 90.15    &  & 91.83  &  & 93.53  &  & 90.77    &  & 92.11  \\
avg.                   &                   & 90.83  &  & 89.74    &  & 90.78  &  & 91.83  &  & 90.18    &  & 91.40  &  & 92.49  &  & 90.55    &  & 91.41  \\ \hline
\end{tabular}
}
\label{T_UTIL}
\end{table}

As depicted in Figure \ref{F_Util 1}, integrating freight into the system enhances passenger utilization, particularly in scenario P+F+M. Interestingly, the overall passenger utilization increases as the peak-to-off-peak ratio increases. Table \ref{T_UTIL} provides further details on average passenger utilization in each period under different scenarios and peak-to-off-peak ratios. Within both peak periods, noticeable increases in passenger utilization indicate that the offered capacity is better utilized as the peak-to-off-peak ratio increases, contributing significantly to the overall enhancement in passenger utilization as depicted in Figure \ref{F_Util 1}.

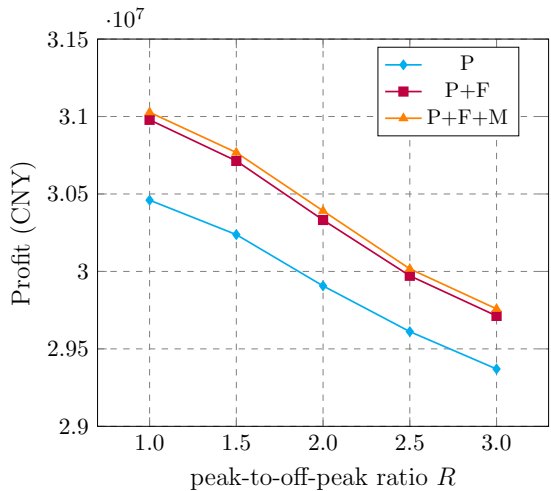
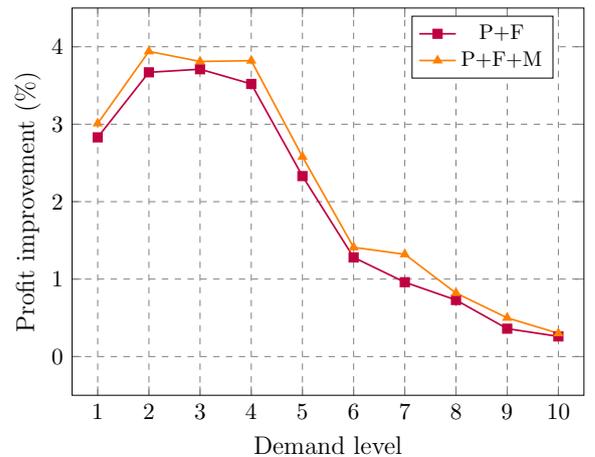
\begin{figure}[!b]
    \centering
    \begin{subfigure}{0.45\textwidth}
        \centering
         \begin{tikzpicture}[scale=0.8]
        \begin{axis}[
            symbolic x coords={1.0,1.5,2.0,2.5,3.0},
            ylabel={Profit (CNY)},
            xlabel={peak-to-off-peak ratio $R$},
            ymin=29000000, ymax=31500000,
            xtick=data,
            enlarge x limits={0.15},
            legend style={font=\footnotesize, xshift=-15mm, anchor=north},
            ylabel style={yshift=-1mm},
            xlabel style={yshift=0mm},
            grid=major,
            grid style={dashed, gray},
            every axis plot/.append style={thick},
            tick label style={font=\small},
            width=9cm,
            height=8cm
        ]
            \addplot[color=cyan, mark=diamond*, mark options={fill=cyan}, thick] coordinates {(1.0, 30459985.15) (1.5, 30237665.32) (2.0, 29906509.44) (2.5, 29611456.22) (3.0, 29370191.34)};
            \addplot[color=purple, mark=square*, mark options={fill=purple}, thick] coordinates {(1.0, 30978630.57) (1.5, 30713430.50) (2.0, 30332048.16) (2.5, 29971873.69) (3.0, 29712938.48)};
            \addplot[color=orange, mark=triangle*, mark options={fill=orange}, thick] coordinates {(1.0, 31027270.93) (1.5, 30767618.54) (2.0, 30390877.38) (2.5, 30016698.88) (3.0, 29758232.23)};
            \legend{P, P+F, P+F+M}
        \end{axis}
    \end{tikzpicture}
    \caption{Average profit for different ratios}
    \label{F_Pro1}
    \end{subfigure}
    \hspace{0.05\textwidth}
    \begin{subfigure}{0.45\textwidth}
        \centering
         \begin{tikzpicture}[scale=0.8]
            \begin{axis}[
                xlabel={Demand level},
                ylabel={Profit improvement (\%)},
                xmin=0.5, xmax=10.5,
                ymin=-0.5, ymax=4.5,
                xtick={1,2,3,4,5,6,7,8,9,10},
                legend style={font=\footnotesize, xshift=-15mm, anchor=north},
                ylabel style={yshift=-1mm},
                xlabel style={yshift=0mm},
                grid=major,
                grid style={dashed, gray},
                tick label style={font=\small},
                width=10cm,
                height=8cm,
                cycle list name=color list
            ]
                \addplot[color=purple, mark=square*, mark options={fill=purple}, thick] coordinates {(1,2.83) (2,3.67) (3,3.71) (4,3.52) (5,2.33) (6,1.28) (7,0.96) (8,0.73) (9,0.36) (10,0.26)};
                \addplot[color=orange, mark=triangle*, mark options={fill=orange}, thick] coordinates {(1,3.01) (2,3.94) (3,3.81) (4,3.82) (5,2.58) (6,1.41) (7,1.32) (8,0.82) (9,0.50) (10,0.30)};
                \legend{P+F, P+F+M}
            \end{axis}
        \end{tikzpicture}
     \caption{Average profit improvement for increasing demand}
    \label{F_Pro2}
    \end{subfigure}
    \caption{Average profit (CNY) / profit improvement (\%) obtained by integrating freight under different scenarios}
    \label{F_Pro}
\end{figure}

Figure \ref{F_Util 2} shows that passenger demand patterns do not significantly impact freight utilization, which remains high. The lowest freight utilization is 96.8\% in scenario P+F and 97.9\% in scenario P+F+M. Moreover, integrating freight into the system is more effective when mixed trains are allowed, as evidenced by the higher freight utilization under scenario P+F+M.

The level of service and the utilization of trains directly impact the system's profitability. Figure \ref{F_Pro1} illustrates the average profits across ten demand levels under different scenarios as a function of the peak-to-off-peak ratio. As anticipated, integrating freight services enhances the profitability of the system. Scenario P+F+M demonstrates yields higher profits compared to scenario P+F, as it allows serving more passengers and freight with higher train loading rates.

Figure \ref{F_Pro2} further shows how the profit improvement, averaged across the different peak-to-off-peak ratios, varies with increasing demand levels. Here, `profit improvement' refers to the increase in system profit achieved by integrating freight compared to the total profit when only passengers are transported. Both curves in Figure \ref{F_Pro2} follow a similar pattern: an initial ascent followed by a descent. Initially, as demand increases, the network can more efficiently pool passengers and freight, resulting in an initial rise in profit improvement. However, with further increases in passenger demand, there is less slack in system capacity, such that the advantage of being able also to transport freight diminishes. 

\subsection{Benefits of implementing a multi-period plan}

Many railway systems are based on a single line plan that is repeatedly executed across all periods, here referred to as the `periodic (PE)' approach, which can be convenient both from an operational standpoint and from the perspective of passengers. However, our study allows for a unique plan for each period, known as the `multi-period (MP)' approach. This section explores the advantages of integrating an MP line plan compared to a PE line plan. The PE line plan is derived by introducing constraints into the MP model to ensure that if a line is operated in one period, it is operated at the same frequency in all other periods.

\begin{figure}[!b]
    \centering
    \begin{subfigure}{0.3\textwidth}
        \centering
         \begin{tikzpicture}[scale=0.8]
         \begin{axis}[
            ybar, 
            scale only axis,
            bar width=9pt,
            height=5cm,
            width=5.5cm,
            xlabel={peak-to-off-peak ratio $R$},
            ylabel={Profit (CNY)},
            ymin=24000000, ymax=32500000,
            xtick=data,
            symbolic x coords={1,1.5,2,2.5,3},
            legend style={font=\tiny, at={(0.98,0.97)}, anchor=north east, legend columns=1},
            ylabel style={yshift=-5pt},
            grid=major,
            grid style={dashed, gray},
        ]
        \addplot[color=red, fill=red!30]
            coordinates {(1,30469363) (1.5,29286976.6199999) (2,28099383.0208333) (2.5,27243262.0999999) (3,26638087.115)};
        \addlegendentry{PE}
        \addplot [color=blue, fill=blue!30]
            coordinates {(1,30459985.145) (1.5,30237665.3183333) (2,29906509.435) (2.5,29611456.2216666) (3,29370191.335)};
        \addlegendentry{MP}
        \end{axis}
    \end{tikzpicture}
    \caption{Scenario P}
    \end{subfigure}
    \hspace{0.03\textwidth}
    \begin{subfigure}{0.3\textwidth}
        \centering
         \begin{tikzpicture}[scale=0.8]
         \begin{axis}[
            ybar, 
            scale only axis,
            bar width=9pt,
            height=5cm,
            width=5.5cm,
            ylabel={Profit (CNY)},
            xlabel={peak-to-off-peak ratio $R$},
            ymin=24000000, ymax=32500000,
            xtick=data,
            symbolic x coords={1,1.5,2,2.5,3},
            legend style={font=\tiny, at={(0.98,0.97)}, anchor=north east, legend columns=1},
            ylabel style={yshift=-5pt},
            grid=major,
            grid style={dashed, gray},
        ]
         \addplot[color=red, fill=red!30]
            coordinates {(1,30580516.5591428) (1.5,29410531.5350819) (2,28258647.43) (2.5,27348366.545) (3,26736190.2916667)};
        \addlegendentry{PE}
        \addplot[color=blue, fill=blue!30]
            coordinates {(1,30978630.57) (1.5,30713430.495) (2,30332048.1583333) (2.5,29971873.6883333) (3,29712938.4775)};
        \addlegendentry{MP}
        \end{axis}
    \end{tikzpicture}
    \caption{Scenario P+F}
    \end{subfigure}
     \hspace{0.03\textwidth}
     \begin{subfigure}{0.3\textwidth}
        \centering
         \begin{tikzpicture}[scale=0.8]
         \begin{axis}[
            ybar, 
            scale only axis,
            bar width=9pt,
            height=5cm,
            width=5.5cm,
            ylabel={Profit (CNY)},
            xlabel={peak-to-off-peak ratio $R$},
            ymin=24000000, ymax=32500000,
            xtick=data,
            symbolic x coords={1,1.5,2,2.5,3},
            legend style={font=\tiny, at={(0.98,0.97)}, anchor=north east, legend columns=1},
            ylabel style={yshift=-5pt},
            grid=major,
            grid style={dashed, gray},
        ]
        \addplot[color=red, fill=red!30]
            coordinates {(1,30686617.795) (1.5,29495628.2172126) (2,28369201.4422225) (2.5,27388766.81) (3,26775946.135)};
        \addlegendentry{PE}
        \addplot[color=blue, fill=blue!30]
            coordinates {(1,31027270.9333333) (1.5,30767618.54) (2,30390877.3758333) (2.5,30016698.88) (3,29758232.2325)};
        \addlegendentry{MP}  
        \end{axis}
    \end{tikzpicture}
    \caption{Scenario P+F+M}
    \end{subfigure}
    \caption{Average profit (CNY) across ten demand levels  under different ratios and scenarios obtained by MP and PE line plans}
    \label{F PROFIT}
\end{figure}
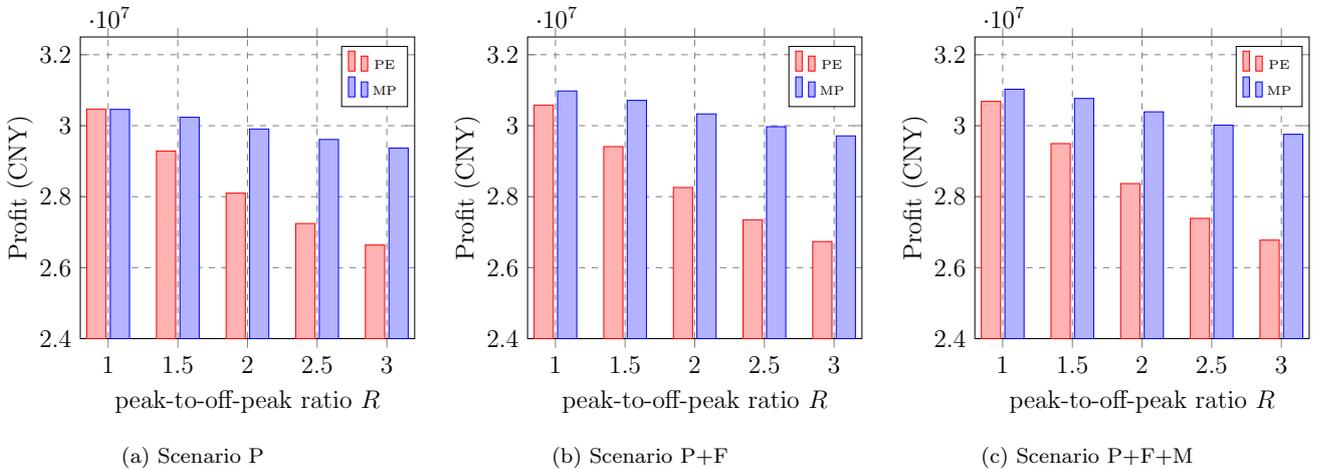

We highlight the profitability benefits of adopting MP line planning. Figure \ref{F PROFIT} presents three subplots that compare the average profits of MP and PE line plans under the different peak-to-off-peak ratios and scenarios. When the demand pattern is uniform, both approaches achieve the same profit when serving only passengers, and minor increases in profit can be observed in MP plans when integrating freight services. As the peak-to-off-peak ratio increases, the benefit of the MP approach becomes increasingly apparent in all scenarios, as this approach provides a tailored line plan for each period. 

\begin{table}[!t]
\renewcommand{\arraystretch}{0.8}
\centering
\caption{Average service levels (\%) across ten demand levels under different ratios and scenarios in MP and PE line plans}
\adjustbox{max width=1\textwidth}{
\begin{tabular}{ccccccccccccccccccccccc}
\hline
\multirow{3}{*}{Demand}      &  & \multirow{3}{*}{Scenario} &  & \multicolumn{19}{c}{peak-to-off-peak ratio of passenger demand $R$}                                                                                              \\ \cline{5-23} 
                           &  &                           &  & \multicolumn{3}{c}{1.0} &  & \multicolumn{3}{c}{1.5} &  & \multicolumn{3}{c}{2.0} &  & \multicolumn{3}{c}{2.5} &  & \multicolumn{3}{c}{3.0} \\ \cline{5-7} \cline{9-11} \cline{13-15} \cline{17-19} \cline{21-23} 
                           &  &                           &  & PE       &    & MP      &  & PE       &    & MP      &  & PE       &    & MP      &  & PE       &    & MP      &  & PE       &    & MP      \\ \hline
\multirow{3}{*}{Passenger} &  & P                         &  & 88.87    &    & 88.85   &  & 87.38    &    & 88.51   &  & 85.56    &    & 87.38   &  & 85.00    &    & 86.86   &  & 84.12    &    & 85.76   \\
                           &  & P+F                       &  & 88.96    &    & 88.91   &  & 87.35    &    & 88.55   &  & 85.51    &    & 87.37   &  & 85.51    &    & 86.89   &  & 84.34    &    & 85.83   \\
                           &  & P+F+M                     &  & 89.18    &    & 89.11   &  & 88.36    &    & 88.79   &  & 88.36    &    & 88.46   &  & 86.40    &    & 87.38   &  & 84.65    &    & 86.14   \\ \hline
\multirow{2}{*}{Freight}   &  & P+F                       &  & 6.24     &    & 24.74   &  & 5.86     &    & 23.46   &  & 7.06     &    & 22.57   &  & 2.25     &    & 20.99   &  & 2.56     &    & 20.81   \\
                           &  & P+F+M                     &  & 8.50     &    & 25.03   &  & 4.11     &    & 24.64   &  & 6.30     &    & 23.08   &  & 3.15     &    & 22.34   &  & 2.82     &    & 22.25   \\ \hline
\end{tabular}
}
\label{T_service}
\end{table}

Table \ref{T_service} compares passenger and freight service levels for MP and PE line plans under different scenarios and peak-to-off-peak ratios. Regarding passenger service levels, MP and PE plans show similar service levels when demand is uniformly distributed across periods. However, as the peak-to-off-peak ratio increases, the difference becomes more pronounced but consistently remains within 2\%, which is still significant given the large base of passenger demand. Regarding freight service levels, the highest service level under the PE plan reaches only 8.5\%, while the lowest under the MP plan reaches 20.8\%. The maximum difference in percentage points is observed in scenario P+F+M with a ratio of 1.5, exceeding 20\%.

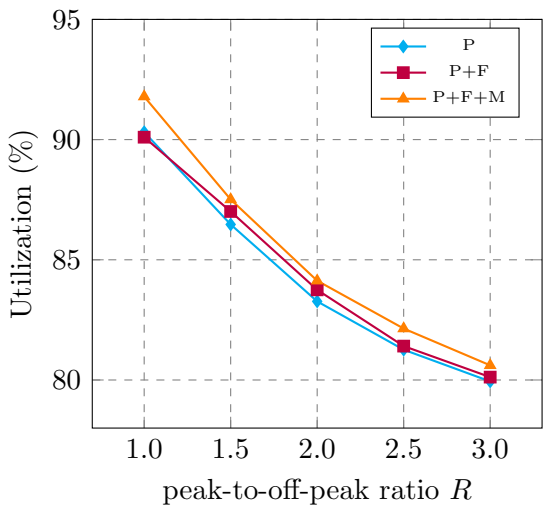
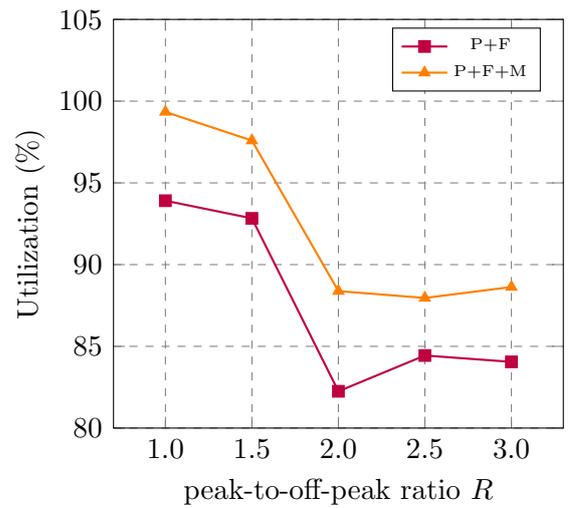
\begin{figure}[htbp]
    \centering
    \begin{subfigure}{0.45\textwidth}
        \centering
         \begin{tikzpicture}[scale=1]
        \begin{axis}[
            xlabel={peak-to-off-peak ratio $R$},
            ylabel={Utilization (\%)},
            symbolic x coords={1.0,1.5,2.0,2.5,3.0},
            ymin=78, ymax=95,
            xtick=data,
            enlarge x limits={0.15},
            legend style={font=\tiny, at={(0.95,0.98)}, anchor=north east},
            ylabel style={font=\small, yshift=-1mm},
            xlabel style={font=\small, yshift=0mm},
            grid=major,
            grid style={dashed, gray},
            every axis plot/.append style={thick},
            tick label style={font=\small},
            width=7.5cm,
            height=7cm
        ]

            \addplot[color=cyan, mark=diamond*, mark options={fill=cyan}] coordinates {
                (1.0,90.31) (1.5,86.47) (2.0,83.27) (2.5,81.26) (3.0,79.94)
            };
            \addplot[color=purple, mark=square*, mark options={fill=purple}] coordinates {
                (1.0,90.1) (1.5,87.01) (2.0,83.75) (2.5,81.41) (3.0,80.12)
            };
            \addplot[color=orange, mark=triangle*, mark options={fill=orange}] coordinates {
                (1.0,91.79) (1.5,87.51) (2.0,84.13) (2.5,82.14) (3.0,80.61)
            };
            \legend{P, P+F, P+F+M}
         \end{axis}
    \end{tikzpicture}
    \caption{Passenger total utilization}
    \label{F_PE UTIL 1}  
    \end{subfigure}
    \hspace{0.05\textwidth}
    \begin{subfigure}{0.45\textwidth}
        \centering
        \begin{tikzpicture}[scale=1]
        \begin{axis}[
            xlabel={peak-to-off-peak ratio $R$},
            ylabel={Utilization (\%)},
            symbolic x coords={1.0,1.5,2.0,2.5,3.0},
            ymin=80, ymax=105,
            xtick=data,
            enlarge x limits={0.15},
            legend style={font=\tiny, at={(0.95,0.98)}, anchor=north east},
            ylabel style={font=\small, yshift=-1mm},
            xlabel style={font=\small, yshift=0mm},
            grid=major,
            grid style={dashed, gray},
            every axis plot/.append style={thick},
            tick label style={font=\small},
            width=7.5cm,
            height=7cm
        ]
        
            \addplot[color=purple, mark=square*, mark options={fill=purple}] coordinates {
                (1.0,93.91) (1.5,92.83) (2.0,82.25) (2.5,84.44) (3.0,84.05)
            };
            \addplot[color=orange, mark=triangle*, mark options={fill=orange}] coordinates {
                (1.0,99.34) (1.5,97.59) (2.0,88.38) (2.5,87.96) (3.0,88.63)
            };
          
            \legend{P+F, P+F+M}
        \end{axis}
    \end{tikzpicture}
    \caption{Freight total utilization}
    \label{F_PE UTIL 2}  
    \end{subfigure}
        \caption{Average total utilization (\%) across ten demand levels under different ratios and scenarios in PE line plans}
    \label{F_PE UTIL}    
    \end{figure}

In addition to service levels, the efficiency of the PE system is also lower. Figures \ref{F_PE UTIL 1} and \ref{F_PE UTIL 2} depict the average passenger and freight utilization across ten demand levels under various scenarios and peak-to-off-peak ratios. Firstly, there is a notable decrease in passenger utilization as the peak-to-off-peak ratio increases, which contrasts with the rising passenger utilization observed under the MP plan in Figure \ref{F_Util 1} as the peak-to-off-peak ratio increases. To clarify this trend, we examine passenger utilization in each period in Table \ref{T_PE P UTIL}. The decline is primarily attributed to a significant drop in passenger utilization during off-peak periods as the peak-to-off-peak ratio increases. As the ratio increases, more lines are needed at higher frequencies during peak periods, whereas the capacity required during the off-peak periods decrease. Therefore, adopting the same line plan in each period will result in train capacity being underutilized during off-peak periods. Secondly, freight utilization shows a significant decline under the PE plan, dropping from 93.9\% to 82.2\% in scenario P+F and from 99.3\% to 87.9\% in scenario P+F+M, as depicted in Figure \ref{F_PE UTIL 2}. In contrast, under the MP plan, the decrease is much smaller, from 98.1\% to 96.8\% in scenario P+F and from 98.6\% to 97.9\% in scenario P+F+M, as depicted in Figure \ref{F_Util 2}. This highlights the superior robustness of the MP plan in terms of freight utilization.

\begin{table}[!t]
 \renewcommand{\arraystretch}{0.8}
\centering
\caption{Average passenger utilization (\%) of each period across ten demand levels under different ratios and scenarios in PE line plans}
\adjustbox{max width=1\textwidth}{
\begin{tabular}{ccccccccccccccccccc}
\hline
\multirow{2}{*}{$R$} &  & \multicolumn{5}{c}{Scenario P}       &  & \multicolumn{5}{c}{Scenario P+F}       &  & \multicolumn{5}{c}{Scenario P+F+M}       \\ \cline{3-7} \cline{9-13} \cline{15-19} 
                       &  & M-peak &  & Off-peak &  & E-peak &  & M-peak &  & Off-peak &  & E-peak &  & M-peak &  & Off-peak &  & E-peak \\ \hline
1.0                      &  & 90.71  &  & 90.14    &  & 90.07  &  & 90.18  &  & 90.02    &  & 90.10  &  & 91.91  &  & 91.92    &  & 91.53  \\
1.5                    &  & 94.30  &  & 70.36    &  & 94.52  &  & 95.17  &  & 70.49    &  & 95.18  &  & 94.93  &  & 72.37    &  & 94.97  \\
2.0                     &  & 94.96  &  & 60.39    &  & 95.16  &  & 94.97  &  & 59.38    &  & 94.75  &  & 95.21  &  & 61.10    &  & 95.45  \\
2.5                    &  & 94.80  &  & 52.72    &  & 95.06  &  & 95.18  &  & 52.95    &  & 95.23  &  & 95.27  &  & 54.71    &  & 95.34  \\
3.0                      &  & 95.11  &  & 48.76    &  & 95.03  &  & 95.27  &  & 49.40    &  & 94.99  &  & 95.51  &  & 50.97    &  & 95.45  \\
avg.                   &  & 93.98  &  & 64.47    &  & 93.97  &  & 94.15  &  & 64.45    &  & 94.05  &  & 94.57  &  & 66.21    &  & 94.55  \\ \hline
\end{tabular}
}
\label{T_PE P UTIL}
\end{table}

\subsection{Impact of demand patterns and modes on the line plan}

This section evaluates the impact of passenger demand patterns and operated modes on the line plan, analyzing changes in the number of lines, total frequency, and proportions of different types of trains.

\begin{figure}[!t]
    \centering
    \begin{subfigure}{0.45\textwidth}
        \centering
         \begin{tikzpicture}[scale=1]
        \begin{axis}[
            xlabel={peak-to-off-peak ratio $R$},
            ylabel={Number of lines},
            symbolic x coords={1.0,1.5,2.0,2.5,3.0},
            ymin=80, ymax=130,
            xtick=data,
            enlarge x limits={0.15},
            legend style={font=\tiny, at={(0.95,0.98)}, anchor=north east},
            ylabel style={font=\small, yshift=-1mm},
            xlabel style={font=\small, yshift=0mm},
            grid=major,
            grid style={dashed, gray},
            every axis plot/.append style={thick},
            tick label style={font=\small},
            width=7.5cm,
            height=7cm
        ]

             \addplot[color=cyan, thick, mark=diamond*, mark options={fill=cyan}]
                coordinates {(1.0,92.6) (1.5,91.6) (2.0,88.9) (2.5,86.4) (3.0,85.9)};
            \addplot[color=purple, thick, mark=square*, mark options={fill=purple}]
                coordinates {(1.0,109.7) (1.5,109.1) (2.0,106.8) (2.5,102.2) (3.0,101)};
            \addplot [color=orange, thick, mark=triangle*, mark options={fill=orange}]
                coordinates {(1.0,120.9) (1.5,118.9) (2.0,112.7) (2.5,109.6) (3.0,109.2)};
            \legend{P, P+F, P+F+M}
        \end{axis}
    \end{tikzpicture}
    \caption{Average number of lines}
    \label{line}  
    \end{subfigure}
    \hspace{0.05\textwidth}
    \begin{subfigure}{0.45\textwidth}
        \centering
        \begin{tikzpicture}[scale=1]
        \begin{axis}[
            xlabel={peak-to-off-peak ratio $R$},
            ylabel={Total frequency},
            symbolic x coords={1.0,1.5,2.0,2.5,3.0},
            ymin=340, ymax=460,
            xtick=data,
            enlarge x limits={0.15},
            legend style={font=\tiny, at={(0.95,0.98)}, anchor=north east},
            ylabel style={font=\small, yshift=-1mm},
            xlabel style={font=\small, yshift=0mm},
            grid=major,
            grid style={dashed, gray},
            every axis plot/.append style={thick},
            tick label style={font=\small},
            width=7.5cm,
            height=7cm
        ]
        \addplot[color=cyan, thick, mark=diamond*, mark options={fill=cyan}]
                coordinates {(1.0,369.1) (1.5,362.7) (2.0,356.8) (2.5,352.9) (3.0,348.5)};
            \addplot [color=purple, thick, mark=square*, mark options={fill=purple}]
                coordinates {(1.0,439.6) (1.5,426.5) (2.0,416.4) (2.5,406.8) (3.0,402.6)};
            \addplot [color=orange, thick, mark=triangle*, mark options={fill=orange}]
                coordinates {(1.0,443) (1.5,431.7) (2.0,421.6) (2.5,410.7) (3.0,407.8)};
            \legend{P, P+F, P+F+M}
        \end{axis}
    \end{tikzpicture}
    \caption{Average total frequency}
    \label{fre}  
    \end{subfigure}
    \caption{Average number of lines and total frequency across ten demand levels under different ratios and scenarios}
    \label{line fre}    
    \end{figure}
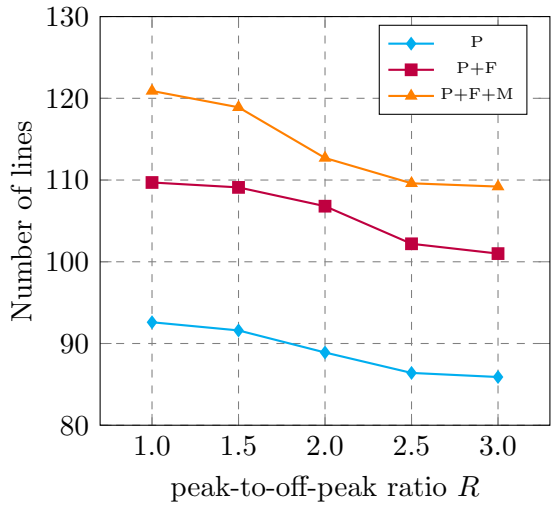
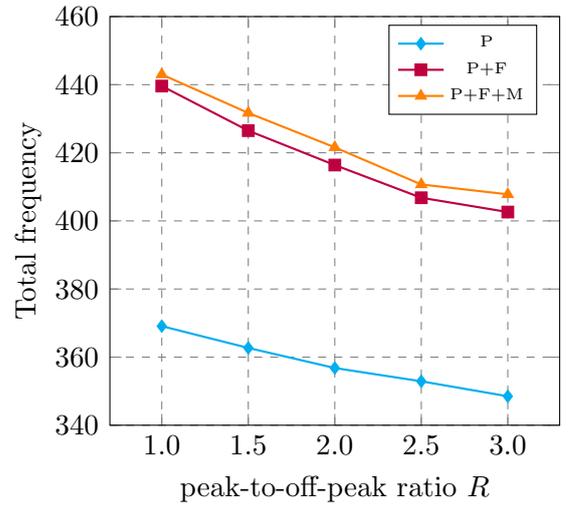

The two subplots in Figure \ref{line fre} illustrate the average numbers of installed lines and total frequencies across ten demand levels under the different peak-to-off-peak ratios and scenarios. More operated modes results in more installed lines, which could increase the complexity of the operations. However, it also boosts the total frequency, signaling improved system service capacity. Additionally, system capacity is highest when the passenger demand is uniform.

\begin{table}[!b]
\renewcommand{\arraystretch}{0.8}
\centering
    \caption{Average proportion (\%) of trains operated in different modes across ten demand levels for different ratios under scenarios P+F and P+F+M}
\adjustbox{max width=1\textwidth}{
\begin{tabular}{ccccccccccc}
\hline
\multirow{2}{*}{$R$} & \multirow{2}{*}{} & \multicolumn{3}{c}{Scenario P+F} &  & \multicolumn{5}{c}{Scenario P+F+M}     \\ \cline{3-5} \cline{7-11} 
                       &                   & Passenger train  &   & Freight train  &  & Passenger train  &  & Mixed train  &  & Freight train \\ \hline
1.0                      &                   & 83.34       &   & 16.66    &  & 74.18  &  & 17.97 &  & 7.85    \\
1.5                    &                   & 84.69       &   & 15.31    &  & 75.52  &  & 16.79 &  & 7.69    \\
2.0                      &                   & 85.40       &   & 14.60    &  & 77.23  &  & 15.56 &  & 7.21    \\
2.5                    &                   & 85.92       &   & 14.08    &  & 77.27  &  & 15.47 &  & 7.26    \\
3.0                      &                   & 86.21       &   & 13.79    &  & 79.87  &  & 13.67 &  & 6.46    \\
Avg.                   &                   & 85.11       &   & 14.89    &  & 76.82  &  & 15.89 &  & 7.29    \\ \hline
\end{tabular}
}
\label{T Lin plan}
\end{table}

Furthermore, Table \ref{T Lin plan} provides the average proportions of different types of trains in Scenarios P+F and P+F+M under the different peak-to-off-peak ratios. The dominance of passenger trains highlights their central role as the primary emphasis of the system's service. Especially as the peak-to-off-peak ratio increases, there is a decrease in the proportion of mixed trains and freight trains and an increase in the proportion of passenger trains, indicating that the integrated passenger-freight railway system increasingly prioritizes passengers over freight.

\section{Conclusions and future research}
This paper proposes a strategic solution for utilizing surplus capacity in railway systems by investigating a multi-period line planning problem for integrated passenger-freight transportation. We develop an MIP model based on a period-extended change\&go-network that allows tracing the lines and transfers taken by passengers and freight, as well as the period in which they travel. We propose two column generation-based heuristics, price-and-branch, and a diving method to solve practical instances.

Computational results based on a Chinese high-speed railway network show that the diving heuristic outperforms the price-and-branch method in terms of solution quality and computational time, achieving average gaps below 1.61\% in medium instances and below 2.84\% in large instances. The integration of freight does considerably increase the complexity of the problem, as the time to solve the freight pricing problem is a large contributor to the overall solving time of the diving algorithm.

The findings of this study suggest that integrating freight into railway passenger operations can enhance the system's profitability and operational efficiency, especially when mixed trains, that combine passengers and freight, are considered. While periodic line plans are perhaps easier to operate and understand for passengers, multi-period line plans can better accommodate the time-varying nature of passenger demand, increasing profitability and efficiency. Introducing additional train modes can lead to more installed lines, increasing operational complexity but also enhancing service capacity through higher frequencies. Across all studied passenger demand patterns, the primary focus of the system remains on serving passengers, as the vast majority of installed lines are operated using dedicated passenger trains.

For future research, it would be interesting to consider the interaction of line planning and timetabling for an integrated passenger-freight rail system and assess to what extent the proposed line planning model facilitates better timetables. In terms of extending the line planning model, one could consider stochastic patterns in passenger and/or freight demand to better capture real-life variability, or an extension towards rolling stock planning to better assess the operational costs of a line plan.

\clearpage
\bibliographystyle{plainnat}
\bibliography{ref}

\clearpage
\appendix 
\section{Details in demand generation}
The stations considered in the network are categorized into four classes based on their scale:

\begin{itemize}
     \item Major stations: Nanjing Nan, Suzhou, Shanghai, Shanghaihongqiao, Hangzhou Dong.
    \item Intermediate stations: Wuxi, Kunshan Nan, Jiaxing Nan.
    \item Small stations: Zhenjiang, Danyang, Changzhou, Suzhouuyuanqu, Changzhou Bei, Wuxi Dong, Suzhou Bei, Jiangning, Wawushan, Yixing, Huzhou.
    \item Minor stations: Stations not mentioned above.
 \end{itemize}

 \begin{table}[ht]
\centering
\caption{Expected total demand (in Carriages) between stations of different classes during the plan horizon}
\adjustbox{max width=1\textwidth}{
\begin{tabular}{ccccccccccc}
\hline
\multicolumn{5}{c}{Medium network}                  &  & \multicolumn{5}{c}{Large network}                       \\ \cline{1-5} \cline{7-11} 
\multicolumn{1}{c|}{Class} & Major   & Intermediate   & Small    & Minor    &  & \multicolumn{1}{c|}{Class} & Major   & Intermediate   & Small   & Minor   \\ \cline{1-5} \cline{7-11} 
\multicolumn{1}{c|}{Major}     & 24   & 20 & 16   & 12  &  & \multicolumn{1}{c|}{Major}     & 24  & 16    & 12  & 8    \\
\multicolumn{1}{c|}{Intermediate }     & 20 & 16   & 12 & 8    &  & \multicolumn{1}{c|}{Intermediate }     & 16   & 12  & 8    & 6 \\
\multicolumn{1}{c|}{Small}     & 16   & 12 & 8   & 4  &  & \multicolumn{1}{c|}{Small}     & 12 & 8    & 6 & 2 \\
\multicolumn{1}{c|}{Minor}     & 12 & 8   & 4 & 2 &  & \multicolumn{1}{c|}{Minor}     & 8   &6 & 2 & 2 \\ \hline
\end{tabular}
}
\label{demand}
\end{table}
 
Table \ref{demand} presents the expected total demand between stations of different classes over the planning period in two networks, measured in carriages. Due to the differing numbers of OD pairs considered in medium and large networks, the expected level of demand varies accordingly. We assume a distribution of 70\% for passenger demand and 30\% for freight demand. For each OD pair, we first identify the respective station classes and obtain the corresponding expected values. Passenger and freight demand is then randomly generated, allowing for fluctuations of ±10\% for passenger demand and ±5\% for freight demand based on these expected values. For example, the expected total demand from Nanjing South to Wuxi is 20 carriages, consisting of 14 carriages for passenger demand and 6 carriages for freight demand. Consequently, the passenger demand will be randomly generated within the range of 12.6 to 15.4 carriages, while the freight demand will range from 5.7 to 6.3 carriages. Finally, the specific quantities of passenger and freight demand will be determined based on the generated number of carriages and the loading capacities per carriage.
\end{document}